\documentclass[%
aip,
amsmath,amssymb,
preprint,%
]{revtex4-1}
\usepackage[dvipsnames]{xcolor}
\usepackage{graphicx}% Include figure files
\usepackage{dcolumn}% Align table columns on decimal point
\usepackage{bm}% bold math
\usepackage{xr}
\usepackage{amssymb}
\usepackage{amsmath}
\usepackage{mathrsfs}

\usepackage{bm}
\usepackage{mathrsfs}
\usepackage{amssymb}
\usepackage{graphicx}
\usepackage{subfigure}
\usepackage{color}
\usepackage{amsmath,graphicx,subfigure,float,appendix,color}
\usepackage{subfigure}
\usepackage{graphicx}
\usepackage{color}
\usepackage{mathrsfs}
\usepackage{appendix}
\DeclareMathOperator*{\esssup}{ess\,sup}

\usepackage[pdfborder={0 0 0},colorlinks=true,citecolor=blue,linkcolor=blue,CJKbookmarks=true,citecolor=blue]{hyperref}
\usepackage[utf8]{inputenc}
\usepackage[T1]{fontenc}
\usepackage{mathptmx}
\usepackage{etoolbox}
\usepackage{subfigure}
\usepackage[ruled,lined]{algorithm2e}
\makeatletter
\def\@email#1#2{%
	\endgroup
	\patchcmd{\titleblock@produce}
	{\frontmatter@RRAPformat}
	{\frontmatter@RRAPformat{\produce@RRAP{*#1\href{mailto:#2}{#2}}}\frontmatter@RRAPformat}
	{}{}
}%
\makeatother

 \usepackage{tabls}
 \newtheorem{theorem}{THEOREM}[section]
\newtheorem{definition}[theorem]{Definition}
\newtheorem{assumption}[theorem]{Assumption}
\newtheorem{remark}[theorem]{Remark}
\newtheorem{corollary}[theorem]{Corollary}
\newtheorem{lemma}[theorem]{Lemma}
\newtheorem{proposition}[theorem]{Proposition}
\newtheorem{example}[theorem]{Example}

\begin{document}
	
	\preprint{AIP/123-QED}
	
\title{Invariance Principles for $G$-Brownian-Motion-Driven Stochastic Differential Equations and Their Applications to $G$-Stochastic Control}
\author{Xiaoxiao Peng}	
\affiliation{School of Mathematical Sciences, Shandong University, Shandong 250100 ,China}
%\affiliation{Research Institute of Intelligent Complex Systems, Fudan University, Shanghai 200433, China}
\author{Shijie Zhou}
\affiliation{Research Institute of Intelligent Complex Systems, Fudan University, Shanghai 200433, China}	
\author{Wei Lin} \email{wlin@fudan.edu.cn}
\affiliation{Research Institute of Intelligent Complex Systems, Fudan University, Shanghai 200433, China}
\affiliation{State Key Laboratory of Medical Neurobiology and MOE Frontiers Center for Brain Science, Institutes of Brain Science, Fudan University, Shanghai 200032, China}
\author{Xuerong Mao}
\affiliation{Department of Mathematics and Statistics, University of Strathclyde, 26 Richmond Street, Glasgow, G11XT, UK}
%\affiliation{State Key Laboratory of Medical Neurobiology and MOE Frontiers Center for Brain Science, Institutes of Brain Science, Fudan University, Shanghai 200032, China}
%\affiliation{School of Mathematical Sciences, LMNS, and SCMS, Fudan University,
%Shanghai 200433, China}

	\date{\today}
	
		\maketitle

%\tableofcontents

\section{abstract}
 The {\it G}-Brownian-motion-driven stochastic differential equations ({\it G}-SDEs) as well as the {\it G}-expectation, which were seminally proposed by Peng and his colleagues, have been extensively applied to describing a particular kind of uncertainty arising in real-world systems modeling.  Mathematically depicting long-time and limit behaviors of the solution produced by {\it G}-SDEs is beneficial to  understanding the mechanisms of system's evolution.  Here, we develop a new {\it G}-semimartingale convergence theorem and further establish a new invariance principle for investigating the long-time behaviors emergent in {\it G}-SDEs.  We also validate the uniqueness and the global existence of the solution of {\it G}-SDEs whose vector fields are only locally Lipschitzian with a linear upper bound.  To demonstrate the broad applicability of our analytically established results, we investigate its application to achieving {\it G}-stochastic control in a few representative dynamical systems.

\section{Introduction}	
Long-time and limit behaviors of the solutions generated by stochastic differential equations (SDEs) have received growing attention because such behaviors usually correspond to particular  functions in real-world systems \cite{gevers1991continuous, Mao-5,Mao-6, florchinger1995lyapunov, brockwell1999stability}. Interesting physical or/and biological phenomena have been systematically investigated, including asymptotic behaviors of random matrices in quantum physics \cite{mendelson2014singular}, stochastic resonance \cite{benzi1983theory}, stochastic homogeneity \cite{DabrockHofmanova-520}, stochastic stabilization or synchronization \cite{b33,b34,li2019robust,li2011output}, and random-temporal-structure-induced emergence \cite{b38,b39,b40,b41}.  Also developed were stochastic versions of invariance principle, which originated from LaSalle's invariance principle \cite{Lasalle-1,Lasalle-2} for deterministic systems and then has been extended successfully to study the SDEs \cite{Mao-7,zhou2023generalized}, the stochastic differential delayed equations (SDDEs) \cite{XuerongMao-8, Mao-9}, the stochastic functional differential equations (SFDEs) \cite{Mao-10, YiQi-11} and even the discrete stochastic dynamical systems \cite{zhou2022generalized}.  These versions of invariance principles are often used to elucidate the asymptotic behaviors, such as stability, boundedness, and invariance in some chaotic attractors, emergent in random systems.

In addition to the traditional frameworks of randomness and stochasticity,  measuring uncertainties of randomness is another important issue in those areas replete with fluctuations and 
risks of high level, such as economics \cite{Knight-12}. A seminal framework by means of sublinear expectation was fundamentally built by Peng and his colleagues to quantify such uncertainties \cite{Peng-13} and then extended broadly in line with the modern probability theory.  Indeed, the framework has been put forward to investigating the $G$-Brownian-motion-driven stochastic differential equations ($G$-SDEs), which thus provides a 
model to describe the randomness with uncertainties in evolutionary dynamics.  Also systematically investigated was the well-posedness of $G$-SDEs \cite{Gao-14,Peng-13} and stochastic functional differential equations ($G$-SFDEs) \cite{YongRen-15,FaizullahFaiz-18}.
Furthermore, although the stability of $G$-SDEs has been widely investigated \cite{LiLin-16,RenYin-17}, rigorously delicate descriptions of stability, boundedness, control and even invariance property in dynamical attractors using $G$-SDEs are still lacking.

In this article, we, therefore, intends to fill in this gap through novelly developing an invariance principle for $G$-SDEs and investigate its applicability to the stochastic control, especially in the case that the noise is uncertain.  As such,  this invariance principle can render the analytical investigations of dynamics produced by $G$-SDEs much clearer and more complete.   In order to develop this new principle, we need to establish a new version of $G$-semimartingale convergence theorem, nontrivially generalizing the classical semimartingale convergence theorem developed in \cite{LiptserShiryayev-24}.  	

The remaining of this article is organized as follows. 
Section~\ref{pre} introduces some basic concepts and provides some preliminary theorems of sublinear expectations. 
Section \ref{Gsemi} rigorously proves the $G$-semimartingale convergence theorem as follows.
\begin{theorem}\label{semimartinagel-intro}
	Assume $A^{1}$ and $A^{2}$ are two non-decreasing process with initial value 0, $A^{1}(t) $ is a continuous process and $\mathbb{\hat{E}}[A^{1}(+\infty)]<+\infty $. Assume that $Z$ is a non-negative $G$-semimartingale satisfying $\mathbb{\hat{E}}[Z^{+}({0})]<\infty$ with the form as
	$
		Z({t})=Z({0})+A^{1}({t})-A^{2}(t)+M({t}), \  t\geq 0,
$
	where $M({t})$ is a continuous $G$-supermartingale with initial value 0. $M({t})\in L_{G}^{1}(\Omega_{t})$ for every $t\geq 0$. Then, we have that $A^{2}(+\infty)<+\infty$, $\lim_{t\rightarrow +\infty}Z({t})$ finitely exists, and that $\lim_{t\rightarrow +\infty}M({t})$ finitely exists quasi-surely.
\end{theorem}	

\noindent Here, we sketch the proof of the above convergence theorem as follows. By extending the space of random variables, we generalize Fatou's Lemma on the $G$-conditional expectation.   Combining with the uppercrossing inequality, we derive the $G$-martingale convergence theorem for a continuous process and then establish the essential $G$-semimartingale convergence theorem. Also in this section, we present the other more applicable versions of the $G$-semimartingale convergence theorem. 	 With all these preparations, Section \ref{secGSDE} presents our main result, the \textit{invariance principle} for the $G$-SDEs, and validates it using the established $G$-semimartingale convergence theorem.   Here, we show this principle as follows.

\begin{theorem}\label{lasalle}
	With those conditions and assumptions listed in Section \ref{secGSDE},  we suppose that there exists a function $V \in C^{2,1}(\mathbb{R}^{d}\times \mathbb{R}_{+}; \mathbb{R}_{+})$, a function $\gamma \in L^{1}(\mathbb{R}_{+}; \mathbb{R}_{+})$ and a continuous function $\eta: \mathbb{R}^{d} \rightarrow \mathbb{R}_{+}$ such that 
	$\lim_{\vert x\vert \rightarrow \infty} \inf_{0\leq t <+\infty} V({\bm x},t)=\infty$
	and 
	$\mathcal{L}V({\bm x},t)\leq \gamma(t)-\eta({\bm x})$,
	where the diffusive operator 
	$
	\mathcal{L}V=V_{t}+V_{x_{i}}f^{i}+G\Big((V_{x_{k}}(h^{kij}+h^{kji})+V_{x_{k}x_{l}} g^{ki}g^{lj})_{i,j=1}^{n}\Big )
	$
	where Einstein's notations are applied here.  Then, we have that
	$\lim_{t \rightarrow +\infty}V({\bm x}(t),t)$ finitely exists quasi-surely
	and that
	$\lim_{t \rightarrow +\infty}\eta({\bm x}(t))=0$  quasi-surely.
	Moreover, we have 
	$
	\lim_{t \rightarrow +\infty}d(\bm x(t), {\rm Ker}(\eta))=0.
	$
	Here, $\bm{x}(t)$ is the solution of the $G$-SDEs which read
	\begin{equation}
		{\rm d}{\bm x}(t)={\bm f}({\bm x}(t),t){\rm d}t+{\bm g}({\bm x}(t),t){\rm d}{\bm B}(t)+{\bm h}({\bm x}(t),t){\rm d}\langle {\bm B}\rangle (t).
	\end{equation}
\end{theorem}

\noindent The proof of such theorem, though inspired by \cite{Mao-7}, is rather different. By $G$-Itô's formula, we write out the function in a form of the $G$-semimartingale and then apply the corresponding convergence theorem.   By estimating the calculus of $\eta$ based on the uppercrossing stopping time, we show that all trajectories converge to the kernel of the function $\eta$ quasi-surely. Still in this section, we further present several generalized versions of invariance principle.   All these build up a solid foundation for Section \ref{example}, where we use the {\it G}-stochastic control to stabilize  representative complex dynamics, demonstrating the broad applicability of our analytically-established results.  Finally, Section \ref{discussion} provides some discussion and concluding remarks.
%	At last, some discussions and conclusions are collected in section \ref{discussion}.

\section{Preliminaries}\label{pre}

In this section, we present some frequently used definitions and results of sublinear expectation theory, which will be useful for our following investigations.   For more details, we refer to \cite{LaurentDenis-25,Peng-13,Peng-26, Lipeng-27}. 	

To begin with, we let $\Omega$ be a given set, and $\mathcal{H}$ be the space of all real-valued functions defined on $\Omega$. Denote by ${C}_{l,{\rm Lip}}(\mathbb{R}^{d })$ the space of all locally {\rm {\rm Lip}}schitz-continuous functions on $\mathbb{R}^{d }$. And, for any function $\varphi \in {C}_{l,{\rm Lip}}(\mathbb{R}^{d })$, if $x_{i}(\omega) \in \mathcal{H}$ for all $i=1,2,\cdots,d$, then $\varphi(x_{1}(\omega), \cdots, x_{d}(\omega))\in \mathcal{H}$.

Next, we provide some basic concepts on the 
sublinear expectation.

\begin{definition}[Sublinear Expectation \cite{Peng-13}]
	A functional $\mathbb{{E}}[\cdot]$ is said to be a sublinear expectation on $\mathcal{H}$ if it satisfies:
		{\rm (1)} $\mathbb{E}[c]=c, \ {\rm for} \ {\rm any}~c \in \mathbb{R}$, {\rm (2)} $\mathbb{E}[X] \leq \mathbb{E}[Y], \ {\rm for} \ {\rm any}~X\leq Y$,  {\rm (3)} $\mathbb{E}[X+Y]\leq \mathbb{E}[X]+\mathbb{E}[Y]$, and {\rm (4)} $\mathbb{E}[\lambda X]=\lambda\mathbb{E}[X], \ {\rm for} \ {\rm any} \lambda \geq 0$.
\end{definition}

\begin{definition}[$G$-Function \cite{Peng-13}]\label{def2}
	A function $G:\mathbb{R}^{d} \times \mathbb{S}^{d} \to \mathbb{R}$ is said to be sublinear and monotone if it satisfies
$(1)~G({\bm p}+\bar{\bm p}, {\bm A}+\bar{{\bm A}}) \leq G({\bm p}, {\bm A})+G(\bar{\bm p}, \bar{\bm A})$, 
		$(2)~G({\bm p}, {\bm A}) \leq G({\bm p}, \bar{\bm A}), \text{if}~{\bm A} \leq \bar{\bm A}$, and  $(3)~G(\lambda {\bm p}, \lambda {\bm A})=\lambda G({\bm p}, {\bm A}),~\forall~\lambda \geq 0$.
 
	Here, $\mathbb{S}^{d}$ denotes the space of $d \times  d$ symmetric matrices. And ${\bm A} \leq \bar{\bm A}$ implies the nonnegativity of the symmetric matrix $\bar{\bm A}-{\bm A}$.
\end{definition}

 In the following, we assume the function $G$ defined in Definition \ref{def2} is independent of the vector $p$. It is worthwhile to mention that, when $d=1$, $G$ is reduced to the form $G(r)=\frac{1}{2}(r^{+}\overline{\sigma}^{2}-r^{-}\underline{\sigma}^{2})$ for some \textcolor{black}{non-negative} $\underline{\sigma} \leq \overline{\sigma}$. Here $r^+$ and $r^-$ correspond to the non-negative and the non-positive parts of $r$, respectively.  Moreover, if a symmetric $G$-Brownian motion satisfies ${\mathbb{\hat{E}}}[{\bm A}{\bm B}(t),{\bm B}(t)]=2G(\bm {A})t$ with $G({\bm A})=\frac{1}{2}{\mathbb{\hat{E}}}[{\bm A}{\bm B}(1),{\bm B}(1)]$, then $G$ is said to be a $G$-function related to  the symmetric $G$-Brownian motion ${\bm B}$. Here, the definition of $G$-Brownian motion, as well as $G$-conditional expectation,  can be found in \cite{Peng-13}. 

 Moreover, it is necessary to introduce some definitions on some spaces of functions and measures. Here, we denote, respectively, by

$\bullet$ $\mathscr{F}_{t}:$ The completion of $\sigma({\bm B}(s):s\leq t)$,

$\bullet$ $\mathscr{B}(\Omega):$ The Borel $\sigma$-algebra on $\Omega$,

$\bullet$ $L^{0}(\Omega)$: The space of all $\mathscr{B}(\Omega)$-measurable functions,

$\bullet$  $L_{G}^{p}(\Omega):$ The completion of the space ${\rm Lip}(\Omega)$ under the norm $\|\cdot\|_{L_{G}^{p}}:=(\mathbb{\hat{E}}[\vert \cdot\vert ^{p}])^{\frac{1}{p}}$, 

$\bullet$  $  {\rm Lip}\left(\Omega_{t}\right)$: $\{\varphi({\bm B}({t_{1}}), {\bm B}({t_{2}})-{\bm B}({t_{1}}), \cdots,$ $ {\bm B}({t_{k}})-{\bm B}({t_{k-1}})):\varphi \in {C}_{l,{\rm Lip}}(\mathbb{R}^{m \times k }), ~0\leq t_{1}<\cdots<t_{k}\leq t\}$,

$\bullet$  $ L_{G}^{p}(\Omega_{t}):$ $L_{G}^{p}(\Omega) \cap {\rm Lip}\left(\Omega_{t}\right)$,

$\bullet$  $ \mathcal{M}: $ The set of all probability measure defined on $\Omega$,

$\bullet$  $E_{Q}[\cdot]$: The expectation under the traditional probability measure $Q$,

$\bullet$  \textcolor{black}{$\mathcal{P}(t, Q)$:=} $\left\{R \in \mathcal{M}: E_{Q}[X]=E_{R}[X], \forall X \in {\rm Lip}\left(\Omega_{t}\right)\right\}$,

$\bullet$  $	\mathcal{Q}:=\{Q\in \mathcal{M}:E_{Q}[X]\leq {\mathbb{\hat{E}}}[X], \forall X \in L_{G}^{1}\left(\Omega\right)\}$, and

$\bullet$ \textcolor{black}{$\mathcal{L}^{0}(\Omega)$:=} \{$X \in L^{0}(\Omega)$: $E_{Q}[X]$ exists for any $Q \in \mathcal{Q}$\}.

From Theorem 1.2.1 in \cite{Peng-26}, it follows that the sublinear expectation satisfies $\mathbb{\hat{E}}[X]=\sup_{Q \in \mathcal{Q}}E_{Q}[X]$ for each $X \in$ Lip($\Omega$).  Thus, the definition of $\mathbb{\hat{E}}[\cdot]  $ can be extended to $\mathcal{L}^0(\Omega)$. In addition, for the $G$-conditional expectation defined above, it can be represented by means of the probability space. 

\begin{theorem}[\cite{HuPeng-29}]\label{G-CE}
	For each $Q \in \mathcal{Q}$ and $X \in L_{G}^{1}(\Omega)$,
	$
		{\mathbb{\hat{E}}}_{t}[X]=\esssup_{R \in \mathcal{P}(t, Q)}{ }^{Q} E_{R}\left[X \mid \mathscr{F}_{t}\right], \ \ Q\text{-}a.s..
	$
	Here, if \ $Y=\esssup_{R \in \mathcal{P}(t, Q)}{ }^{Q} E_{R}\left[X \mid \mathscr{F}_{t}\right]$, it means that for every $R \in \mathcal{P}(t, Q)$, $E_{R}\left[X \mid \mathscr{F}_{t}\right] \leq Y, Q\text{-}a.s..$ Moreover, if $E_{R}\left[X \mid \mathscr{F}_{t}\right] \leq Z$ \textcolor{black}{for each $R \in \mathcal{P}(t, Q)$}, $Q\text{-}a.s.$, then we must have $Y\leq Z, Q\text{-}a.s..$
\end{theorem}

For introducing $G$-Itô's calculus, we define $M_{G}^{p}([0, T])$, a space of random process, and the $G$-Itô's calculus on it (refer to \cite{Peng-13} for details). Moreover, the quadratic variation is defined in the same manner as that in normal stochastic analysis.  However, the range of the quadratic variation here is much different.

\begin{lemma}[\cite{Peng-13}]\label{lemmalemma}
	For an $m$-dimensional $G$-Brownian motion ${\bm B}$, there exists a bounded, convex and closed set $\Gamma \in \mathbb{S}_{+}^{m}$ such that 
	$\langle {\bm B} \rangle({t}) \in t\Gamma:=\{t \gamma: \gamma \in \Gamma\}$,
	where ${\mathbb{S}_{+}^{n}}$ represents the space of all positive symmetric matrices. Also, $	\langle {\bm B} \rangle({t})$ and $	\langle {\bm B} \rangle({t+s})-	\langle {\bm B} \rangle({s})$ are identically distributed. 
\end{lemma}

	\begin{remark}\label{lemmaqua}
	In what follows, denote by $\bar{\gamma}:=
	\max_{\gamma \in \Gamma} (\vert \gamma\vert _{F}\vee \vert \gamma\vert_{2} ) $ where $\vert \cdot\vert _{F}$ and $\vert \cdot\vert_{2} $, respectively, are the Frobenius norm \cite{belitskii2013matrix} and 2-norm for the matrix. Then, it follows from Lemma \ref{lemmalemma} that $\vert 	\langle {\bm B} \rangle({t})\vert_{F} \vee \vert \langle {\bm B} \rangle({t})\vert_{2} \leq \bar{\gamma}t$. Especially when $m=1$, we have $\bar{\gamma}=\overline{\sigma}^{2}$. Also, the largest eigenvalue of a matrix is denoted by $\lambda_{\rm max}(\cdot)$. 
\end{remark}

There are some very useful inequalities for our investigation in this article. Combining the results of Sections 3.3-3.5 in \cite{Peng-13}, Lemma \ref{lemmalemma}, and Remark \ref{lemmaqua}, we give the conclusions as follows.

\begin{theorem} \label{BDG}
	For any $\eta(t), ~\gamma(t)\in M_{G}^{2}[0, T]$,  we have
	$
		\mathbb{\hat{E}}\left(\int_{0}^{T} \eta(t) \mathrm{~d} B_{i}(t)\right)=0$ and $\mathbb{\hat{E}}\left(\int_{0}^{T} \eta(t) \mathrm{~d} B_{i}(t)\int_{0}^{T} \gamma(t) \mathrm{~d} B_{j}(t)\right)
		=\mathbb{\hat{E}}\left(\int_{0}^{T} \eta(t)\gamma(t) \mathrm{~d}\langle B_{i},B_{j}\rangle(t)\right)
		\leq \bar{\gamma} \cdot \mathbb{\hat{E}}\left(\int_{0}^{T} \vert \eta(t)\gamma(t)\vert  \mathrm{~d} t\right).
	$
\end{theorem}

Now we introduce the Choquet capacity and some related propositions.
\begin{definition}[Choquet Capacity, \cite{Peng-13}]\label{defcho}
	For $\mathcal{A} \in \mathscr{B}(\Omega)$,  define by
	$
		c(\mathcal{A}):=\sup_{Q \in \mathcal{Q}}Q[\mathcal{A}]=\mathbb{\hat{E}}[1_{\mathcal{A}}].
	$ A property is called \textit{valid quasi-surely} if this property is valid on the set $\Omega \backslash \mathcal{A} $ with $c(\mathcal{A})=0$.
\end{definition}

\begin{proposition}[Monotone Convergence Theorem, \cite{LaurentDenis-25, Peng-26}]\label{monotone}
	If $X({n}) \uparrow X$, $\{X(n)\}\subset \mathcal{L}^{0}(\Omega)$, $X({n})$ is nonnegative, then $\mathbb{\hat{E}}[X({n})]\uparrow \mathbb{\hat{E}}[X]$.
\end{proposition}

\begin{theorem}[\cite{Li-31}]\label{MarCon}
	Assume that $\{M(n) \}$ is a $G$-supermartingale, satisfying $\sup_{n}\mathbb{\hat{E}}[M^{-}(n)]<+\infty$. Then, $\lim_{n\rightarrow \infty}M(n)$ exists, which is finite quasi-surely. Here, the definition of $G$-martingale can be found in \cite{Peng-13}.
\end{theorem}

	\section{$G$-Semimartingale Convergence Theorem}\label{Gsemi}

In the literature, the semimartingale convergence theorem mainly describes the asymptotic property of the semimartingale, which is a random variable comprising a martingale and a process with bounded variation. Inspired by this well-established and broadly-applied convergence theorem, we are to establish a $G$-semimartingale convergence theorem and its variant.  It will be shown that the $G$-semimartingale convergence theorem is based crucially 
on Doob's $G$-martingale convergence theorem.  In fact, to our best knowledge, the continuous version of Doob's $G$-martingale convergence theorem has not yet been established until the result presented \textcolor{black}{as follows}.

\begin{proposition}[$G$-Martingale Convergence Theorem, A Continuous Version]\label{convergence}
	Assume that $\{M({t}):t\in[0,+\infty)\}$ is a right- or left-continuous 
	$G$-supermartingale, and $M(t)\in L_{G}^{1}(\Omega_{t})$. Moreover, assume that $\mathbb{\hat{E}}[\sup_{t\geq 0}M^{-}({t})]<+\infty$.  Then, $M({t})$ converges finitely to $M(+{\infty})\in L_{G}^{1^{*}_{*}}(\Omega)$ quasi-surely. Moreover,  $\mathbb{\hat{E}}_{t}[M({+\infty})]\leq M({t})$. \textcolor{black}{Here, the definition of $L_{G}^{1^{*}_{*}}(\Omega)$ is provided in Definition \ref{remark1} of Appendix~\ref{appendix1}}.
\end{proposition}

\textcolor{black}{The proof of this proposition is tedious and tangential to the main focus of this article. To enhance the readability, we include the proof into Appendix \ref{appendix1}.}   Now, with this preparation, we establish the following $G$-semimartingale convergence theorem.

\begin{theorem}[$G$-Semimartingale Convergence Theorem]\label{semimartinagel}
	Assume that $A^{1}$ and $A^{2}$ are two non-decreasing processes with initial value $0$, and that $A^{1}(t) $ is a continuous process with $\mathbb{\hat{E}}[A^{1}(+\infty)]<+\infty$.  Also, assume that $Z$ is a non-negative $G$-semimartingale satisfying $\mathbb{\hat{E}}[Z^{+}({0})]<\infty$ with the form
	$
		Z({t})=Z({0})+A^{1}({t})-A^{2}(t)+M({t}), \  t\geq 0,
	$
	where $M({t})$ is a continuous $G$-supermartingale with initial value $0$ and $M({t})\in L_{G}^{1}(\Omega_{t})$ for every $t\geq 0$. Then, we have that $A^{2}(+\infty)<+\infty$, $\lim_{t\rightarrow +\infty}Z({t})$ finitely exists and $\lim_{t\rightarrow +\infty}M({t})$ finitely exists quasi-surely.
\end{theorem}

\noindent{\bf Proof.}
Notice that
$
	M({t})=Z({t})-Z({0})-A^{1}(t)+A^{2}(t)\geq -Z({0})-A^{1}(+\infty).
$
Then, $\sup_{t\geq 0}M^{-}({t}) \leq Z^{+}(0)+A^{1}(+\infty)$. By Proposition \ref{convergence}, we have $\lim_{t \rightarrow \infty}M({t})$ finitely exists quasi-surely. Because
$
	A^2({t})=Z({0})+A^{1}({t})+M({t})-Z({t})\leq Z({0})+A^{1}({t})+M({t}) 
$
and 
$
	Z({t})=Z({0})+A^{1}({t})-A^{2}({t})+M({t}),
$
their limits also exist quasi-surely. 

It is mentioned that this $G$-semimartingale convergence theorem can only deal with the case where the limit of $A^{1}({t})$ is supposed to be finite under the sublinear expectation.    We now give its variant, the $G$-semimartingale convergence theorem with the $\mathbb{F}$-stopping time.  It can deal with the case where the condition on the finite limit of $A^{1}({t})$ in Theorem \ref{semimartinagel} is removed.  The tradeoff however requires
more conditions for the $G$-martingale $M$.	

\begin{theorem}[$G$-Semimartingale Convergence Theorem with Stopping Time]\label{semimartinage2}
	Assume that $A^{1}$ and $A^{2}$ are two non-decreasing processes both with initial value $0$, and that $A^{1}({t}) $ is a continuous adapted process.  Also assume that $Z$ is a non-negative adapted process satisfying $\mathbb{\hat{E}}[\vert Z({0})\vert ]<\infty$ with the form
	$
		Z({t})=Z({0})+A^{1}({t})-A^{2}({t})+M({t}), \ \  t\geq 0,
	$
	where $M({t})$ is a continuous process with initial value $0$.  Furthermore, assume that there exists a series of $\mathbb{F}$-stopping times $\tau_{N}$ satisfying $\{\tau_{N}\rightarrow +\infty\}$ quasi-surely such that, for any $Q \in \mathcal{Q}$, ${E}_{Q}[M({t\wedge\tau_{N}})\vert \mathscr{F}_{s}]=M({s\wedge\tau_{N}})$. Then, we have quasi-surely
$$
	\begin{array}{l}
		\displaystyle\left\{\omega:A^{1}(+\infty)<+\infty\right\}\subset  \left\{\omega:\lim_{t\rightarrow +\infty}Z({t})  \mbox{ finitely  exists} \right\} \\
		\displaystyle ~~~~\cap \left\{\omega:A^{2}(+\infty)<+\infty\right\}\cap
		\left\{\omega:\lim_{t\rightarrow +\infty}M({t}) \ \mbox{finitely  exists}\right\}.
	\end{array}
$$
Here, $\mathcal{A} \subset \mathcal{B}$ quasi-surely means that $c(\mathcal{A}\backslash \mathcal{B})=0$, where $c$ is the Choquet capacity provided in Definition \ref{defcho}.
\end{theorem}	

\noindent{\bf Proof.}
Denote by
$\mathcal{A}=\Omega \backslash \Big(\{\omega:\lim_{t\rightarrow +\infty}Z({t})  \  \mbox{finitely  exists} \}\cap \{\omega:A^{2}(+\infty)<+\infty\}\cap\{\omega:\lim_{t\rightarrow +\infty}M({t}) \ \mbox{finitely  exists} \} \Big)$.
For every $Q \in \mathcal{Q}$, we have $E_{Q}[\vert Z({0})\vert ]\leq \mathbb{\hat{E}}[\vert Z({0})\vert ]$. By the $G$-semimartingale convergence theorem for the normal probability space \cite{ LiptserShiryayev-24}, we have $Q(\mathcal{A})=0$. By the arbitrariness of the $Q$'s choice, we obtain that $c(\mathcal{A})=\sup_{Q \in \mathcal{Q}}Q(\mathcal{A})=0$, which therefore completes the proof.

\section{Invariance Principle in Sublinear Expectation}\label{secGSDE}

Now, we consider a $d$-dimensional $G$-stochastic differential equation which reads
\begin{equation}\label{GSDE}
	{\rm d}{\bm x}(t)={\bm f}({\bm x}(t),t){\rm d}t+{\bm g}({\bm x}(t),t){\rm d}{\bm B}(t)+{\bm h}({\bm x}(t),t){\rm d}\langle {\bm B}\rangle (t),
\end{equation}
where the initial value $x(0)=x_{0}$.  Furthermore, we denote,  respectively, by 
$
\vert \bm A\vert _{2}:=\sqrt{{\rm tr}(\bm A^{\top}\bm A)}~~\mbox{ and} ~~\vert \bm A\vert :=\vert \bm A\vert_{F}=\sqrt{\sum_{i,j=1}^{n}a_{ij}^{2}} 
$
different norms of a given matrix $\bm{A}$. 	All  functions ${\bm f}:\mathbb{R}^{d}\times \mathbb{R}_{+}\rightarrow \mathbb{R}^{d}$, ${\bm g}:\mathbb{R}^{d}\times \mathbb{R}_{+} \rightarrow \mathbb{R}^{d\times m}$, and ${\bm h}:\mathbb{R}^{d}\times \mathbb{R}_{+} \rightarrow \mathbb{R}^{d\times m \times m}$ are supposed to be continuous.  In addition,  $h^{kij}=h^{kji}$, and $f^{i}({\bm x},\cdot)$, $g^{ij}({\bm x},\cdot)$ and $h^{kij}({\bm x},\cdot) \in M_{G}^{2}[0,T]$ for every $T>0$. We need the following assumptions.

\begin{assumption}\label{asump1}
	For any $N \in \mathbb{N}$, there exists a number 
	$C_{N}$ such that
	$
		\vert {\bm f}({\bm x},t)-f({\bm y},t)\vert +\vert {\bm g}({\bm x},t)-{\bm g}({\bm y},t)\vert +\vert {\bm h}({\bm x},t)-{\bm h}({\bm y},t)\vert \leq C_{N}\vert {\bm x}-{\bm y}\vert 
	$
	for all $\vert {\bm x}\vert \wedge \vert {\bm y}\vert \leq N$. Here, $\vert {\bm h}\vert$ still represents the norm for $\bm h$ of $d\times m\times m $ dimensions.
\end{assumption} 

\begin{assumption}\label{asump2}
	There exists a number $C_{l}$ such that		$
		\vert {\bm f}({\bm x},t)\vert +\vert {\bm g}({\bm x},t)\vert +\vert {\bm h}({\bm x},t)\vert \leq C_{l}(1+\vert {\bm x}\vert ), 
	$
	for all $({\bm x},t)\in \mathbb{R}^{d}\times \mathbb{R}_{+}$.
	
\end{assumption}

\textcolor{black}{Underlying these assumptions as prerequisites, the solutions of Eq.~\eqref{GSDE} are well-posed from a certain perspective as follows.}

\begin{proposition}\label{exist}
	If Assumption \ref{asump1} holds, there is a global unique solution in a quasi-sure sense on $[0,\tau_{\infty})$, where 
$
		\tau_{\infty}=\lim_{n\rightarrow +\infty}\tau_{N}, \ \tau_{N}:=\inf \{t\geq 0:\vert {\bm x(t)}\vert \geq N\}.
	$For given $N>0$, there exists $\bm x^{N}\in M_{G}^{2}[0,T]$ with $T>0$ such that $\bm x=\bm x^{N}$ on $[0,\tau_{N})$.  Additionally, for ${\bm A}=(a^{ij}):\mathbb{R}^{d}\times \mathbb{R}_{+} \rightarrow \mathbb{R}^{d\times m}$ with $a^{ij}({\bm x},\cdot)\in M_{G}^{1}[0,T]$ and $T>0$, we have ${\bm M}(t)=\int_{0}^{t\wedge \tau_{N}}{\bm A}({\bm x}(s),s) {\rm d}{\bm B}(s)$ is $Q$-martingale for each $Q \in \mathcal{Q}$. If Assumption \ref{asump2} holds, we have $\tau_{\infty}=+\infty$ quasi-surely.
\end{proposition}

\begin{remark}\label{rk4.4}
\textcolor{black}{ The proof of Proposition \ref{exist} is similar to those presented in Refs.~\cite{Mao-6,LiLin-16}, which we omit here.  It is worth mentioning that $\bm x(\cdot)$, the solution to Eq.~\eqref{GSDE}, does not belong to $M_{G}^{2}([0,T];\mathbb{R}^{d})$. Actually, $\bm x(\cdot\wedge \tau_N)\in M_{*}^{2}([0,T];\mathbb{R}^{d})$ for each $N>0$, which implies that our solution is locally integrable. In particular, if $\tau_\infty =+\infty$, we have $\bm x(\cdot) \in M_{w}^{2}([0,T];\mathbb{R}^{d})$ and it is globally integrable on $[0,+\infty)$ now. Here, both $M_{*}^{2}([0,T];\mathbb{R}^{d})$ and $M_{w}^{2}([0,T];\mathbb{R}^{d})$ are expanded integrand space defined in Chapter 8 of Ref.~\cite{Peng-13} satisfying $M_{G}^{2}([0,T];\mathbb{R}^{d})\subset M_{*}^{2}([0,T];\mathbb{R}^{d})\subset M_{w}^{2}([0,T];\mathbb{R}^{d})$. }
\end{remark}

Next, we introduce $G$-Itô's formula which is useful in the following discussions.
\begin{theorem}[$G$-Itô's formula \cite{Lipeng-27}]\label{GIto}
	Let $V \in C^{2,1}(\mathbb{R}^{d}\times \mathbb{R}_{+}; \mathbb{R}_{+})$. 	For the $d$-dimensional $G$-stochastic differential equations
	$
	{\rm d}{\bm x}(t)={\bm f}(t){\rm d}t+{\bm g}(t){\rm d}{\bm B}(t)+{\bm h}(t){\rm d}\langle {\bm B}\rangle (t)
	$
	with the initial value ${\bm x}( 0)={\bm x}_{0}$. Moreover, ${\bm f}:\mathbb{R}_{+}\rightarrow \mathbb{R}^{d}$, ${\bm g}: \mathbb{R}_{+} \rightarrow \mathbb{R}^{d\times m}$, and ${\bm h}:\mathbb{R}_{+} \rightarrow \mathbb{R}^{d\times m^2}$ with $f^{i}(\cdot),~g^{ij}(\cdot) \in M_{G}^{1}[0,T]$, $h^{kij}(\cdot) \in M_{G}^{2}[0,T]$ for every $T>0$. Then,
	$
	V({\bm x}(t),t)=V({\bm x}_{0},0)+\int_{0}^{t} V_{t}({\bm x}(s),s){\rm d}s+\int_{0}^{t} V_{x_{i}}({\bm x}(s),s)f^{i}(s){\rm d}s
	+ \int_{0}^{t} V_{x_{i}}({\bm x}(s),s)g^{ij}(s){\rm d}B_{j}(s)
	+ \int_{0}^{t} V_{x_{k}}({\bm x}(s),s)h^{kij}(s){\rm d}\langle B_{i}, B_{j} \rangle(s)
	+\int_{0}^{t}\frac{1}{2} V_{x_{k}x_{l}}({\bm x}(s),s)g^{ki}(s)g^{lj}(s){\rm d}\langle B_{i}, B_{j} \rangle(s).
	$
\end{theorem}

\textcolor{black}{Actually,  $G$-Itô's formula presented above could be applicable to $M_{*}^{2}([0,T];\mathbb{R}^{d})$ and $M_{w}^{2}([0,T];\mathbb{R}^{d})$ according to Theorem 5.4 established in \cite{Lipeng-27}. By virture of $G$-Itô's formula,}   Assumption \ref{asump2} used above can be replaced.  To present this result, we introduce  the notation as
$
	\mathcal{L}V:=V_{t}+V_{x_{i}}f^{i}+G\Big((V_{x_{k}}(h^{kij}+h^{kji})+V_{x_{k}x_{l}} g^{ki}g^{lj})_{i,j=1}^{n}\Big),
$
where the function $V \in C^{2,1}(\mathbb{R}^{d}\times \mathbb{R}_{+}; \mathbb{R}_{+})$.    As such, we obtain the following result.

\begin{proposition}\label{exist2}
	Suppose that Assumption \ref{asump1} holds and that there exists a function $\gamma \in L^{1}(\mathbb{R}_{+}; \mathbb{R}_{+})$ such that 
	$
		\mathcal{L}V({\bm x},t)\leq \gamma(t).
	$
	Moreover, $V$ satisfies 
	\begin{equation}\label{conditionofV2}
		\lim_{\vert x\vert \rightarrow \infty} \inf_{0\leq t <+\infty} V({\bm x},t)=+\infty.
	\end{equation}
	Then, $\tau_{\infty}$, as defined in Proposition \ref{exist}, satisfies $\tau_{\infty}=+\infty$ quasi-surely.
\end{proposition}

For simplicity of expression, we still include the proof of Proposition \ref{exist2} in Appendix \ref{appendix3}, where the following proposition is needed.

\begin{proposition}[\cite{LiLin-16}]\label{prop5}
	Let $
	M(t)=\int_{0}^{t} \kappa_{ij}(s){\rm d}\langle B_{i},B_{j} \rangle (s)-\int_{0}^{t} 2G(\bm \kappa){\rm d}s,
	$
	where  ${\bm \kappa} \in M_{G}^{1}([0,T];\mathbb{S}^{n})$.  Then, we have $M(t)\leq 0$ quasi-surely. Particularly $\mathbb{\hat{E}}[M(t)]\leq 0$.
\end{proposition}

  \textcolor{black}{In addition,	we present the following $G$-stochastic Barbalat's lemma that  will be used later, and its proof is provided in Appendix~\ref{appendix_barbalat}}.

  \begin{lemma} \label{theorembarbalat}
  \textcolor{black}{Suppose that Assumption \ref{asump1} holds and $\tau_{\infty}=+\infty $ quasi-surely, where $\tau_\infty$ is defined in Proposition \ref{exist}.  Also suppose that the solution to Eq.~\eqref{GSDE} satisfies $\sup_{t\in \mathbb{R}^{+}}\vert \bm x(t) \vert<+\infty \ q.s.$. Besides, there exists $\eta\in C(\mathbb{R}^{d}; \mathbb{R}_{+})$ such that 
\begin{equation}\label{integralbound}
	\int_{0}^{+\infty}\eta(\bm x(t)){\rm d}t<+\infty, \ \ q.s..
\end{equation}
Then, we have $\lim_{t\rightarrow +\infty}\eta(\bm x(t))=0$ quasi-surely.}
  \end{lemma}

Now, with the following assumption, we state our main theorem.

\begin{assumption}\label{asump3}
	For each $N>0$, $t \in \mathbb{R}_{+}$ and all $\vert \bm x\vert \leq N$, there exists a number $K_{N}>0$ such that 
	$
		\vert {\bm f}({\bm x},t)\vert +\vert {\bm g}({\bm x},t)\vert +\vert {\bm h}({\bm x},t)\vert \leq K_{N}.
	$
\end{assumption}

\begin{theorem}\label{lasalle1}
	Suppose that Assumptions \ref{asump1} and \ref{asump3} hold.   Also suppose that there exist three functions $V \in C^{2,1}(\mathbb{R}^{d}\times \mathbb{R}_{+}; \mathbb{R}_{+})$,  $\gamma \in L^{1}(\mathbb{R}_{+}; \mathbb{R}_{+})$ and $\eta\in C(\mathbb{R}^{d}; \mathbb{R}_{+})$ such that $({\rm UB})$
$\lim_{\vert x\vert \rightarrow \infty} \inf_{0\leq t <+\infty} V({\bm x},t)=\infty$ and 
	$\mathcal{L}V({\bm x},t)\leq \gamma(t)-\eta({\bm x})$.
	Then, we have that $\lim_{t \rightarrow +\infty}V({\bm x}(t),t)$ finitely exists quasi-surely and that
	\begin{equation}\label{limitofeta}
		\lim_{t \rightarrow +\infty}\eta({\bm x}(t))=0 \ \ q.s..
	\end{equation}
	Moreover, $\lim_{t \rightarrow +\infty}d(\bm x, {\rm Ker}(\eta))=0$, where $d(\bm x, {\rm Ker}(\eta)):= \inf_{{\bm y}\in {\rm Ker}(\eta)}\vert \bm x- \bm y \vert $.
\end{theorem}

\noindent{\bf Proof.}
Using Proposition \ref{exist2}, the $G$-SDEs satisfying the conditions assumed in this theorem have a global solution on $[0,+\infty)$ with a property that
$
	\mathcal{L}V({\bm x},t)\leq \gamma(t)-\eta({\bm x})\leq \gamma(t).
$
By $G$-Itô's formula in Theorem \ref{GIto}, Proposition \ref{exist} and Remark \ref{rk4.4}, 
we have
\begin{eqnarray*}
	&&V({\bm x}(t\wedge \tau_{N}),t\wedge \tau_{N})=V({\bm x}_{0},0)+\int_{0}^{t\wedge \tau_{N}} V_{t}({\bm x}(s),s){\rm d}s\\
	&&+\int_{0}^{t\wedge \tau_{N}} V_{x_{i}}({\bm x}(s),s)f^{i}({\bm x}(s),s){\rm d}s  
	+ \int_{0}^{t\wedge \tau_{N}} V_{x_{i}}({\bm x}(s),s)g^{ij}({\bm x}(s),s){\rm d}B_{j}(s)\\
	&&+ \int_{0}^{t\wedge \tau_{N}} V_{x_{k}}({\bm x}(s),s)h^{kij}({\bm x}(s),s){\rm d}\langle B_{i}, B_{j} \rangle(s)
	+\int_{0}^{t\wedge \tau_{N}}\frac{1}{2} V_{x_{k}x_{l}}({\bm x}(s),s)g^{ki}({\bm x}(s),s)\\
	&&g^{lj}({\bm x}(s),s){\rm d}\langle B_{i}, B_{j} \rangle(s),
\end{eqnarray*}
where $\tau_{N}:=\inf \{t\geq 0:\vert {\bm x(t)}\vert \geq N\}$. Letting $N \rightarrow+\infty$ and setting 
$
\bm \kappa=(\kappa_{ij})_{i,j=1}^{m} 
$
for every $t\geq 0$ where
$
	{\kappa}_{ij}=V_{x_{k}}(h^{kij}+h^{kji})+V_{x_{k}x_{l}}g^{ki}g^{lj},
$
we get that $\tau_{N}$ tends to $+\infty$ by Proposition \ref{exist} and
\begin{eqnarray*}
	V({\bm x}(t),t)&=&V({\bm x}_{0},0)+\int_{0}^{t} V_{t}({\bm x}(s),s){\rm d}s+\int_{0}^{t} V_{x_{i}}({\bm x}(s),s)f^{i}({\bm x}(s),s){\rm d}s  \\
	&&+ \int_{0}^{t} V_{x_{i}}({\bm x}(s),s)g^{ij}({\bm x}(s),s){\rm d}B_{j}(s)+\int_{0}^{t} \frac{1}{2}\kappa_{ij}({\bm x}(s),s) {\rm d}\langle B_{i},B_{j} \rangle(s). 
\end{eqnarray*}
Thus, if we set 
$$
	V({\bm x}(t),t)=V({\bm x}_{0},0) + \int_{0}^{t} \gamma(s){\rm d}s-A_{2}(t) +\int_{0}^{t} V_{x_{i}}({\bm x}(s),s)g^{ij}({\bm x}(s),s){\rm d}B_{j}(s),
$$
then  $A_{2}(0)=0$.   Besides,  according to Proposition \ref{prop5},  for every $0\leq t_{1}<t_{2}<+\infty$, we have
\begin{eqnarray*}
&&	A_{2}(t_{2})-A_{2}(t_{1})=\int_{t_{1}}^{t_{2}} \gamma(s){\rm d}s-\int_{t_{1}}^{t_{2}} V_{x_{i}}({\bm x}(s),s)f^{i}({\bm x}(s),s){\rm d}s\\ && -\int_{t_{1}}^{t_{2}} V_{t}({\bm x}(s),s){\rm d}s-\int_{t_{1}}^{t_{2}} \frac{1}{2}\kappa_{ij}({\bm x}(s),s) {\rm d}\langle B_{i},B_{j} \rangle (s)\\
&\geq& \int_{t_{1}}^{t_{2}} \gamma(s){\rm d}s-\int_{t_{1}}^{t_{2}} V_{t}({\bm x}(s),s){\rm d}s\\
	&&-\int_{t_{1}}^{t_{2}} V_{x_{i}}({\bm x}(s),s)f^{i}({\bm x}(s),s){\rm d}s-\int_{t_{1}}^{t_{2}} G(\eta ({\bm x}(s),s)) {\rm d}s\\
	&=&\int_{t_{1}}^{t_{2}} \gamma(s){\rm d}s- \int_{t_{1}}^{t_{2}} \mathcal{L} V({\bm x}(s),s){\rm d}s
	\geq\int_{t_{1}}^{t_{2}} \eta(s){\rm d}s \geq 0
\end{eqnarray*}
which implies that $A_{2}(t)$ is a non-decreasing process.  Using Proposition \ref{exist}, we obtain that $\int_{0}^{t\wedge \tau_{N}} V_{x_{i}}({\bm x}(s),s)g^{ij}({\bm x}(s),s){\rm d}B_{j}(s)$ is a $Q$-martingale for every $Q \in \mathcal{Q}$.   Noticing $\int_{0}^{+ \infty} \gamma(s){\rm d}s<+\infty$ and according to Proposition \ref{semimartinage2}, we have a set ${\Omega_{0}}\subset \Omega$ such that $c({\Omega}\backslash\Omega_{0})=0$. Then, we have that, for all $\omega \in \Omega_{0}$, 
$
	\lim_{n \rightarrow +\infty}A_{2}(t) \ {\rm finitely} \  {\rm exists} 
$
and
$\lim_{n \rightarrow +\infty}V({\bm x}(t),t)$ finitely exists.
Thus, on ${\Omega_0}$, $
	\int_{0}^{+\infty}\eta(\bm x(t)){\rm d}t<+\infty.$  From the finite existence of the limit of $V$, we obtain that, on $\Omega_0$,
$\sup_{t \geq 0}V({\bm x}(t;\omega),t)<+\infty$.
Hence, from the above-assumed condition (UB), it follows that there exists $K(\omega)$ such that
$\sup_{t\geq 0}\vert {\bm x}(t;\omega)\vert \leq K(\omega)$. \textcolor{black}{According to Lemma \ref{theorembarbalat}, we obtain $\lim_{t\rightarrow +\infty} \eta(\bm x(t))=0$ quasi-surely.}

For every $\omega$ satisfying  $\lim_{t\rightarrow +\infty}\eta(\bm x(t;\omega))=0$ and $\sup_{t \in \mathbb{R}_{+}}\vert \bm x(t;\omega)\vert <+\infty$, there exists $\bm y(\omega)$ and a sequence $\{t_{i}\}$ having
$\lim_{i\rightarrow +\infty}\bm x(t_{i};\omega)=\bm y(\omega)$.
So,
$
	\lim_{i\rightarrow +\infty}\eta(\bm x(t_{i};\omega))=\bm \eta(\bm y(\omega))=0
$
and ${\rm Ker}(\eta)\neq \emptyset$. If
$
	\limsup_{t\rightarrow +\infty} d(\bm x(t;\omega),{\rm ker}(\eta))
$ is positive,
there exist a sequence $\{t_{i}\}$ such that
$
	d(\bm x(t_{i};\omega),{\rm ker}(\eta))\geq \epsilon,
$
for some $\epsilon>0$. This implies $\eta(\bm y)>0$,
which is a contradiction.

\begin{remark}
	\textcolor{black}{Here, our conclusions nontrivially extend the corresponding results obtained for the traditional SDEs.  Particularly, the significant differences do exist. First, in terms of the conclusions, we are able to induce relevant results even when the system randomness itself is uncertain, greatly surpassing the applicability scope of existing Brownian motion-driven stochastic systems. From a technical standpoint, our generalized stochastic differential equation (i.e., G-SDE) cannot measure the occurrence probability of events from the perspective of traditional probability measures, but the capacities instead.  Second, the construction of the monotone functions in our semi-martingales differs significantly from the invariance principles in the traditional stochastic analysis. }
\end{remark}

Next, we present another version of invariance principle, where $\eta$ is a function with respect to the function $V$.

\begin{theorem}\label{lasalle2}
	Suppose that Assumption \ref{asump1} holds, and that there exist three functions $V \in C^{2,1}(\mathbb{R}^{d}\times \mathbb{R}_{+}; \mathbb{R}_{+})$, $\gamma \in L^{1}(\mathbb{R}_{+}; \mathbb{R}_{+})$ and $\eta\in C(\mathbb{R}_{+} ;\mathbb{R}_{+})$ such that 
	$
		\mathcal{L}V({\bm x},t)\leq \gamma(t)-\eta(V({\bm x},t)) 
	$
	for all $({\bm x},t)\in \mathbb{R}^{d}\times \mathbb{R}_{+}$. Then, we obtain that $\lim_{t \rightarrow +\infty}V({\bm x}(t),t)$ finitely exists quasi-surely	and 
	$
	\lim_{t \rightarrow +\infty}\eta(V({\bm x}(t),t))=0 \ q.s..
	$	
	Moreover, $\lim_{t \rightarrow +\infty}d( V(\bm x(t),t), {\rm Ker}(\eta))=0.
	$
\end{theorem}

\noindent{\bf Proof.}  Analogously, the $G$-SDEs have a global solution on $[0,+\infty)$ according to Proposition \ref{exist2}.  By the arguments akin to those for validating Theorem \ref{lasalle1}, we obtain
$
	V({\bm x}(t),t)	= V({\bm x}_{0},0) + \int_{0}^{t} \gamma(s){\rm d}s-A_{2}(t)
	+\int_{0}^{t} V_{x_{i}}({\bm x}(s),s)g^{ij}({\bm x}(s),s){\rm d}B_{j}(s),
$
where $A_{2}(0)=0$ and for every $0\leq t_{1}<t_{2}<+\infty$, 
\begin{eqnarray*}
	&&A_{2}(t_{2})-A_{2}(t_{1})=\int_{t_{1}}^{t_{2}} \gamma(s){\rm d}s-\int_{t_{1}}^{t_{2}} V_{x_{i}}({\bm x}(s),s)f^{i}({\bm x}(s),s){\rm d}s\\ &&-\int_{t_{1}}^{t_{2}} V_{t}({\bm x}(s),s){\rm d}s-\int_{t_{1}}^{t_{2}} \frac{1}{2}\kappa_{ij}({\bm x}(s),s) {\rm d}\langle B_{i},B_{j} \rangle(s) \\
\end{eqnarray*}
\begin{eqnarray*}
	&\geq& \int_{t_{1}}^{t_{2}} \gamma(s){\rm d}s-\int_{t_{1}}^{t_{2}} V_{x_{i}}({\bm x}(s),s)f^{i}({\bm x}(s),s){\rm d}s\\
	&&-\int_{t_{1}}^{t_{2}} V_{t}({\bm x}(s),s){\rm d}s-\int_{t_{1}}^{t_{2}} G(\eta ({\bm x}(s),s)) {\rm d}s\\
	&=&\int_{t_{1}}^{t_{2}} \gamma(s){\rm d}s- \int_{t_{1}}^{t_{2}} \mathcal{L} V({\bm x}(s),s){\rm d}s
	\geq\int_{t_{1}}^{t_{2}} \eta(V({\bm x}(s),s)){\rm d}s \geq 0.
\end{eqnarray*}
Hence, by the $G$-semimartingale Convergence Theorem \ref{semimartinage2}, there exists $\bar{\Omega}\subset \Omega$ such that $c({\Omega}\backslash \bar{\Omega})=0$.  Furthermore, we have that, on $\bar{\Omega}$,
$$
	\int_{0}^{\infty} \eta(V({\bm x}(t),t)) {\rm d}t<+\infty~~\mbox{and}~~
	\lim_{n \rightarrow +\infty}V({\bm x}(t),t) \ {\rm finitely} \  {\rm exists} .
$$
Now, we claim that, for every $\omega\in \bar{\Omega}$, we have $\lim_{t \rightarrow +\infty}\eta(V({\bm x}(t;\omega),t))=0$.  We validate the claim by contradiction.  If this is not the case, then we have a sequence $\{t_{k}\}$ with $t_{k+1}-t_{k}>1$ and $\epsilon>0$, such that $\eta(V({\bm x}(t_{k};\omega),t_{k}))>\epsilon$. 
Assume $\sup_{t\geq 0}V({\bm x}(t;\omega),t)\leq K(\omega)$. Hence, there exists $\delta_{1}$ such that
$\vert \eta(x)-\eta(y)\vert \leq \frac{\epsilon}{2}$ for  $0\leq x,y\leq K(\omega)$ and $\vert x-y\vert \leq\delta_{1}$.
As  $\lim_{t\rightarrow +\infty}V({\bm x}(t;\omega),t)$ finitely exists and $V({\bm x}(t;\omega),t)$ is continuous about $t$, we can easily check that it is uniformly continuous on $\mathbb{R}^+$. Thus, there exists $\delta_{2}<1$ such that 
$\vert V({\bm x}(t;\omega),t)-V({\bm x}(s;\omega),s)\vert < \delta_{1}, \ \ \vert t-s\vert < \delta_{2}.
$
Consequently, for $t_{k}\leq t < t_{k}+\delta_{2}$, we have 
$
		\eta(V({\bm x}(t;\omega),t))\geq \eta(V({\bm x}(t_{k};\omega),t_{k}))-\vert \eta(V({\bm x}(t_{k};\omega),t_{k}))
	-\eta(V({\bm x}(t;\omega),t))\vert\geq \frac{\epsilon}{2}.
$
Therefore,
$
	+\infty>\int_{0}^{\infty} \eta(V({\bm x}(t),t)) {\rm d}t\geq \sum_{k=1}^{+\infty}\int_{t_k}^{t_k+\delta_{2}} \eta(V({\bm x}(t),t)) {\rm d}t
	\geq \sum_{k=1}^{+\infty} \frac{\epsilon \delta_{2}}{2}=+\infty,
$
which indicates a contradiction. Finally, the arguments for proving $\lim_{t \rightarrow +\infty}d( V(\bm x(t),t), {\rm Ker}(\eta))=0$ are the same as those for validating the last conclusion in Theorem \ref{lasalle1}.

\begin{remark}
	A set $\mathcal{A} \in \mathscr{B}(\Omega)$ is said to be invariant if 	
	$
		c\big(\{\exists t\geq 0, ~x(t;{\bm x}_{0})\notin \mathcal{A} \}\big)=0,
	$
	for every ${\bm x}_{0} \in \mathcal{A}$.  Actually, if we suppose some conditions to be valid only in the invariant set $\mathcal{A}$ for Theorems $\ref{lasalle1}$ and $\ref{lasalle2}$,
 the conclusions in these theorems still sustain.  
\end{remark}

Finally, we present two corollaries which can be obtained directly form the invariance principles established above.  These results are related to the stability or the exponential stability of the solution $\bm x(t)$.

\begin{corollary}\label{cor4.8}
	Let Assumption \ref{asump1} hold. Assume further that there exists a function $V \in C^{2,1}(\mathbb{R}^{d}\times \mathbb{R}_{+}; \mathbb{R}_{+})$ such that
	\begin{equation}\label{cor1}
		\mu_{1}(\vert {\bm x}\vert )\leq V({\bm x},t)\leq \mu_{2}(\vert {\bm x}\vert ), \ \ \mathcal{L}V({\bm x},t)\leq -\mu_{3}(\vert {\bm x}\vert ),
	\end{equation}
	where $\mu_{1}$, $\mu_{2}$ and $\mu_{3}$ are three strictly increasing functions in $[0,+\infty)$ with the initial value $0$ and
	$\mu_{1}(r), \mu_{2}(r) \rightarrow +\infty$ as $r \rightarrow +\infty$.	
	Then, we have 
	$\lim_{t \rightarrow +\infty}\vert {\bm x}(t)\vert =0 ~q.s.$.
\end{corollary}	

\noindent{\bf Proof.}  From the condition assumed in \eqref{cor1}, it follows that $\mu_{2}^{-1}(V({\bm x},t))\leq \vert {\bm x}\vert$, which implies
$
	\mathcal{L}V({\bm x},t)\leq -\mu_{3}(\mu_{2}^{-1}(V({\bm x},t))).
$
According to Theorem \ref{lasalle2}, we have
$\lim_{t \rightarrow \infty}\mu_{3}(\mu_{2}^{-1}(V({\bm x}(t),t)))=0 ~ q.s.$, which implies 
$\lim_{t \rightarrow \infty}V({\bm x}(t),t)=0 ~ q.s.$.
Therefore, we have 
$\lim_{t \rightarrow \infty}\mu_{1}(\vert {\bm x}(t)\vert )=0 ~ q.s.$, which finally gives
$\lim_{t \rightarrow \infty} \vert {\bm x}(t)\vert =0 ~ q.s.$.

\begin{corollary}\label{cor2}
	Let Assumption \ref{asump1} hold. Assume further 
	that there exist two functions: $V \in C^{2,1}(\mathbb{R}^{d}\times \mathbb{R}_{+}; \mathbb{R}_{+})$ and $\gamma \in L^{1}(\mathbb{R}_{+};\mathbb{R}_{+} )$, such that		
$
		{\rm e}^{\lambda t}\vert {\bm x}\vert ^{p}\leq V({\bm x}(t),t) \ \ \ {\rm and} \ \ \ \mathcal{L}V({\bm x},t)\leq \gamma(t),
	$
	where $\lambda$ and $p$ are positive numbers.
	Then, we have
$
		\varlimsup_{t \rightarrow +\infty}\frac{1}{t} {\rm log}\vert {\bm x}(t)\vert \leq -\frac{\lambda}{p} \ \ \ q.s..
	$
\end{corollary}

\noindent{\bf Proof.} Set $\eta=0$ in Theorem \ref{lasalle2}. Then,
$\lim_{t \rightarrow +\infty}V({\bm x}(t),t)$ finitely exists quasi-surely.
Further use the condition that ${\rm e}^{\lambda t}\vert {\bm x}\vert ^{p}\leq V({\bm x}(t),t)$.   The proof is therefore complete.

	\section{Illustrative Examples: Applying $G$-invariance principle to achieving $G$-stochastic control}\label{example}

In this section, we use several representative examples to illustrate the applicability of our analytical results to realizing $G$-stochastic control of the unstable dynamical systems.

\begin{example}\label{example1}
	Consider a linear (complex network) system  ${\rm d}{\bm x}(t)={\bm A}{\bm x}(t){\rm d}t$. Here, $\bm{A}=[11,5,2;5,11,2;2,2,14]$.  Then, it is easy to check that ${\lambda}_{\rm max}({\bm A})=18$ and the system is unstable. Now, for a $G$-Brownian motion where $\underline{\sigma}^2=3.5$ and $\overline{\sigma}^2=4$, we choose ${\bm D}={\bm I}_{3}$ and $\bm C = [-19,11,2;11,-19,2;2,2,-10]$ to $G$-stochastically control the linear system as 
 $
		{\bm x}(t)={\bm x}_{0}+\int_{0}^{t} {\bm A}{\bm x}(s){\rm d}s + \int_{0}^{t} {\bm D}{\bm x}(s) {\rm d}{\bm B}(s)+\int_{0}^{t} {\bm C}{\bm x}(s) {\rm d}\langle {\bm B}\rangle(s).
	$
  Choosing $V({\bm x}):=\vert {\bm x}\vert ^2$ yields:
	$
		\mathcal{L}V({\bm x})=2{\bm x}^{\top}{\bm A}{\bm x}+G(2{\bm x}^{\top}{\bm D}^{\top}{\bm D}{\bm x}+4{\bm x}^{\top}{\bm C}{\bm x}).
	$
As ${\lambda}_{\rm max}({\bm C})=-6$, we easily derive that $\mathcal{L}V({\bm x})\leq -2.5\vert {\bm x}\vert ^2$. This, according to  Corollary \ref{cor4.8}, ensures the asymptotic stability of the controlled system in a quasi-sure sense.

Moreover, if we set $V({\bm x},t)={\rm e}^{\lambda t}\vert {\bm x}\vert ^2$, we obtain that
	$
		\mathcal{L}V({\bm x},t)=	
		\mathcal{L}V({\bm x})=\left[{\bm x}^{\top}(2{\bm A}+\lambda {\bm I}_{d}){\bm x}+G(2{\bm x}^{\top}{\bm D}^{\top}{\bm D}{\bm x}+4{\bm x}^{\top}{\bm C}{\bm x})\right]{\rm e}^{\lambda t},
	$
	which, using the parameters $\underline{\sigma}^2=3.5$ and $\overline{\sigma}^2=4$, yields $\mathcal{L}V({\bm x},t)\leq (\lambda-1.5)\vert {\bm x}\vert ^{2}$. If we set $\lambda \leq 1.5$, using Corollary \ref{cor2} gives 
	$\varlimsup_{t \rightarrow +\infty}\frac{1}{t} {\rm log}\vert {\bm x}(t)\vert \leq -0.75 \ q.s.$.
	This clearly illustrates the exponential stability of the controlled system. 
\end{example}

\begin{example}\label{exstability}
	Consider an autonomous system, which reads
${\rm d}{\bm x}(t)={\bm f}({\bm x}(t)){\rm d}t$.
Here, ${\bm f}$ satisfies Assumption \ref{asump1} and $\bm f(\bm 0)=\bm 0$. Moreover, $\bm f$ satisfies one-sided Lipschitz condition, i.e., there exists a number $L>0$ such that  
$
		\langle \bm x,\bm f(\bm x) \rangle \leq L  \vert \bm x  \vert^{2}.
$
There are many systems, not globally Lipschitzian, only satisfying this one-sided Lipschitz condition. For instance, both $f(x)=x-x^3$ and the Lorenz system with
	$\bm{f}(\bm x)=[
		\sigma x_{2}-\sigma x_{1},
		\rho x_{1}-x_{3}x_{1}-x_{2},
		x_{1}x_{2}-\beta x_{3} 
		]^{\top}$
satisfy the one-sided Lipschitz condition.  Now, we apply the $G$-stochastic control to the original dynamics, which yields 
$
		{\rm d}{\bm x}(t)={\bm f}({\bm x}(t)){\rm d}t+k\sum_{j=1}^{m}{\bm x}(t){\rm d}{B}_{j}(t)
$
with $k>\left(-{L}/{c_{-1}}\right)^{{1/2}}$ with $c_{-1}:=G\left((-1)_{i,j=1}^{m}\right)$. Here, $(-1)_{i,j=1}^{m}$ corresponds to an $m \times m$ matrix with all elements are $-1$. Then, the controlled system becomes stochastically stable, whose proof is included in Appendix \ref{app4}.  Take the three-dimensional Lorenz system for example.   We are able to use a one-dimensional $G$-Brownian motion to render the controlled system stable quasi-surely, if we set $m=1$, $c_{-1}=G(-1)=-\frac{1}{2}\underline{\sigma}^{2}$, $L\leq \frac{1}{2}(\sigma +\rho)$, and $k>(\sigma +\rho)^{1/2}\underline{\sigma}^{-1}$.
\end{example}

\begin{example}\label{example3}
Consider an oscillating system 
${\rm d}{\bm x}(t)=\bm C \bm f(\bm x(t)){\rm d}t$,
where $\bm C = [1,1,4;5,-1,4;8,1,0]$ and $\bm f(\bm x)=[-x_{1},\arctan (x_{2}), \tanh (x_{3})]^{\top}$.
Now, we consider the $G$-stochastically controlled system as 
$
{\rm d}{\bm x}(t)=\bm C \bm f(\bm x(t)){\rm d}t+ \bm g(\bm x(t)){\rm d}{\bm B}(t)$,
where ${\bm B}$ is a two-dimensional, independent and identically distributed $G$-Brownian motion with $\bar{\sigma}^{2}=50$ and $\underline{\sigma}^{2}=40$, and $\bm g(\bm x)=[\bm A_{1}\bm x, \bm A_{2}\bm x]$ in which ${\bm A_{1}}=[1,0.5,0;0,1,0;0,0,1]$ and ${\bm A_{2}}=[1,0,0;0,1,0.5;0,0,1]$.  Additionally, the $G$-function of $\bm B$ satisfies $\bm G(\bm M)=\sum_{j=1}^{2}G_{j}(a_{jj})$, where $\bm M=(m_{ij})_{i,j=1}^{2}$ is a two-dimensional matrix, and $G_{j}$ is the $G$-function related to the one-dimensional $G$-Brownian motion $B_{j}$. Set $V(\bm x)=\vert \bm x \vert^{\alpha}$ for some $\alpha>0$. By Appendix~\ref{appendix2}, $\mathbb{R}^{3}\backslash \{\bm 0\}$ is an invariant set of the system. It follows that, on $\mathbb{R}^{3}\backslash \{\bm 0\}$,
\begin{eqnarray*}
\mathcal{L}V(\bm x)    &=&\alpha \vert \bm x\vert^{\alpha-2}\big[-x_{1}^{2}+x_{1}\arctan (x_{2})+4x_{1}\tanh (x_{3})-5x_{1}x_{2}-x_{2}\arctan (x_{2})\\
    & &+4x_{2}\tanh(x_{3})-8x_{3}x_{1}+x_{3}\arctan (x_{2})\big]\\
    & & +\alpha \vert \bm x\vert^{\alpha-4} G\left(\vert \bm x\vert^{2}\bm g^{\top}\bm g+(\alpha-2)\bm g^{\top}\bm x \bm x^{\top}\bm g\right)\\
    &\leq  &\alpha \vert \bm x\vert^{\alpha-2}(-x_{1}^{2}+6\vert x_{1}x_{2}\vert + 12\vert x_{1}x_{3}\vert + 5\vert x_{2}x_{3}\vert) \\
    &&+\sum_{j=1}^{2}\alpha \vert \bm x\vert^{\alpha-4}G_{j}\left(\vert \bm x\vert^{2}\vert \bm A_{j}\bm x\vert^{2}+ (\alpha-2)(\bm x^{\top} \bm A_{j}\bm x)^{2}\right).
\end{eqnarray*}
Notice that $(\bm x^{\top}\bm A_{j}\bm x)^{2}\geq \frac{1}{2} \vert \bm x \vert^{2}\vert \bm A_{j}\bm x \vert^{2}+\frac{1}{8} \vert \bm x \vert^{4}$ and $\bm x^{\top}\bm A_{j}\bm x \leq \frac{5}{4} \vert \bm x \vert^{2}$ for $j=1,2,$ and set $\alpha = \frac{2}{25}$. Then, we obtain 
$
\mathcal{L}V(\bm x)\leq \frac{17}{25} \vert \bm x\vert^{\frac{2}{25}}+\sum_{j=1}^{2}\frac{2}{25}\vert \bm x\vert^{-\frac{98}{25}}G_{j}\Big(\frac{2}{25}(\bm x^{\top} \bm A_{j}\bm x)^{2}-\frac{1}{4} \vert \bm x\vert^{4}\Big )\leq -\frac{3}{25}\vert \bm x\vert^{\frac{2}{25}}
 $. Setting $\eta$  in Theorem \ref{lasalle1} as $\eta(\bm x)=\frac{3}{25}\vert \bm x\vert^{\frac{2}{25}}$ guarantees the quasi-sure stability of the above controlled system.
\end{example}

In Appendix \ref{appendix_numerical}, we further provide a few numerical evidences for illustrating the above examples.  It is emphasized that those numerically-presented results do not represent all the exact solution produced by the $G$-SDEs, but only provide some evidences partially supporting the analytical results obtained in the above examples.  The numerical scheme used there is not complete, so it awaits further development for rigorously approximating the solution of $G$-SDEs.

\section{Conclusion}\label{discussion}

In this article, we have developed several invariance principles for the stochastic differential equations driven by the $G$-Brownian motions. Our work is basically inspired by the seminal works from two directions: one is from the stability theory of the traditional SDEs \cite{Mao-7} and the other is from the fundamentally-innovative works on the sublinear expectation \cite{Peng-13}.  
Our contributions include not only the establishment of the $G$-semimartingale convergence theorem and its variants for the sublinear expectation, but also the establishment of several invariance principles and their applications in investigating the long-term behaviors of $G$-SDEs.  Indeed, we anticipate that our analytical results can be beneficial to understanding and solving the problems associated with uncertain randomness in dynamical systems.

As for the future research directions, the assumption on the linear growth and the locally Lipschitz conditions can be further weakened through restricting the discussion for the operator $\mathcal{L}$ in some specific space.  Also, further development of the invariance principles for the $G$-SDDEs and the $G$-SFDEs could be promoted.  More practically, complete scheme for rigorously approximating the solution produced by the $G$-SFDEs deserves deep investigation.

\section{Appendix}\label{appendix}
\subsection{Proof of Proposition \ref{convergence}}\label{appendix1}
\textcolor{black}{First, we establish Fatou's lemma for the $G$-conditional expectation, which is a prerequisite for our proposition to be demonstrated}.

\begin{lemma}[Fatou's Lemma for $G$-conditional Expectation]\label{lemma2.2}
	$\{X(n)\}\in L_{G}^{1}(\Omega)$ are a series of random vectors, and there exists a random variable $M$ such that $\mathbb{\hat{E}}[\vert M\vert ]<+\infty$ and $X(n)\geq M$ for any $n>0$. Then, 
	$
	\mathbb{\hat{E}}_{t}\left[\varliminf_{n\rightarrow \infty}X(n)\right]\leq 
	\varliminf_{n\rightarrow \infty}\mathbb{\hat{E}}_{t}[X(n)].
	$
\end{lemma}

In order to present the proof for this lemma, we need to extend the space of random variables and make some necessary preparations.

\begin{definition}[\cite{HuPeng-29}]\label{remark1}
	Introduce some extended spaces of random variables as follows:
	
	$
	\begin{array}{l}
		\mathcal{L}_{G}^{{1}^{*}}(\Omega):=\Big\{X\in L^{0}(\Omega):  \exists X({n})\in {L}_{G}^{1}(\Omega)  \ {\rm such} \  {\rm that} \ X({n}) \downarrow X  \Big\}, \\
		{L}_{G}^{{1}^{*}}(\Omega):=\Big\{X\in L^{0}(\Omega): \mathbb{\hat{E}}[\vert  X\vert ]<+\infty, \ \ X \in  \mathcal{L}_{G}^{{1}^{*}}(\Omega) \Big\},\\
		\mathcal{L}_{G}^{{1}^{*}_{*}}(\Omega):= \Big\{X\in L^{0}(\Omega):  \exists X({n})\in L_{G}^{1^{*}}(\Omega)  \ {\rm such} \  {\rm that} \ X({n}) \uparrow X  \Big\},\\
		{L}_{G}^{{1}^{*}_{*}}(\Omega):=\Big\{X\in L^{0}(\Omega): \mathbb{\hat{E}}[ \vert  X\vert ]<+\infty, \ \ X \in \mathcal{L}_{G}^{{1}^{*}_{*}}(\Omega) \Big\}.
	\end{array}
	$
	
	Then, we extend the $G$-conditional expectation on $\mathcal{L}_{G}^{{1}^{*}_{*}} (\Omega)$.  Directly, we have ${L}_{G}^{{1}^{*}}(\Omega) \subset	\mathcal{L}_{G}^{{1}^{*}}(\Omega) \subset \mathcal{L}_{G}^{{1}^{*}_{*}}(\Omega)$ and ${L}_{G}^{{1}^{*}}(\Omega) \subset {L}_{G}^{{1}^{*}_{*}}(\Omega) \subset \mathcal{L}_{G}^{{1}^{*}_{*}}(\Omega)$.
\end{definition}

\begin{lemma}[\cite{HuPeng-29}]\label{lemma2.1}
	Suppose that $\{X(n)\}\subset L_{G}^{1^{*}_{*}}(\Omega)$ is a series of non-decreasing random variables.  Denote by $X:=\lim_{n \to \infty}X(n)$.  Then, we have quasi-surely
	$
	\lim_{n \to \infty}\mathbb{\hat{E}}_{t}[X(n)]=\mathbb{\hat{E}}_{t}[X].
	$
\end{lemma}

\begin{lemma}\label{smalllemma1}
	If $X, Y \in L_{G}^{1}(\Omega)$, then $X\wedge Y \in L_{G}^{1}(\Omega)$ $($resp. $X\vee Y \in L_{G}^{1}(\Omega))$.
\end{lemma}

\noindent{\bf Proof.}
As $X, Y \in L_{G}^{1}(\Omega)$, there exists $\{X_{n}\}$ and $\{Y_{n}\}$ contained in ${\rm Lip}(\Omega)$ such that $\mathbb{\hat{E}}[\vert X(n)-X\vert ]\rightarrow 0$ and $\mathbb{\hat{E}}[\vert Y(n)-Y\vert ]\rightarrow 0$. For $\varphi$, $\psi \in C_{l,{\rm Lip}}(\Omega)$, we have 
$
\varphi \wedge \psi=\frac{\varphi+\psi-\vert \varphi-\psi\vert }{2}\in C_{l,{\rm Lip}}(\Omega).
$
Thus, $X(n)\wedge Y(n) \in {\rm Lip}(\Omega)$. So we derive 
$
\mathbb{\hat{E}}[\vert X\wedge Y-X(n)\wedge Y(n)\vert ]\leq \mathbb{\hat{E}}[\vert X-X(n)\vert ]+\mathbb{\hat{E}}[\vert Y-Y(n)\vert ]\rightarrow 0,
$
which implies $X\wedge Y \in L_{G}^{1}(\Omega)$. The case that $X\vee Y \in L_{G}^{1}(\Omega)$ is analogous.

\begin{lemma}\label{lemma3}
	If $X(n)\in L_{G}^{1}(\Omega)$ and $X(n)$ converges to $X$, and there exists a random variable $M$ such that $\mathbb{\hat{E}}[\vert M\vert ]<+\infty$ and $X(n)\geq M$ for any $n>0$. Then, $X\in \mathcal{L}_{G}^{{1}^{*}_{*}}(\Omega).$
\end{lemma}
\noindent{\bf Proof.}
For any $m,n>0$, by Lemma \ref{smalllemma1}, we obtain that $\inf_{n \leq k\leq m} X(k) \in L_{G}^{1}(\Omega)$. Then, from Definition  \ref{remark1}, it follows that $\inf_{k\geq n} X(k) \in \mathcal{L}_{G}^{{1}^{*}}(\Omega)$.  Also, by the fact that $M\leq\inf_{k\geq n} X(k)\leq X(n)$, we have 
$
\left\vert \inf_{k\geq n} X(k)\right\vert \leq \vert X(n)\vert +\vert M\vert . 
$
Thus, $\mathbb{\hat{E}}[\vert \inf_{k\geq n} X(k)\vert ]\leq +\infty$ and $\inf_{k\geq n} X(k) \in{L}_{G}^{{1}^{*}}(\Omega)$ using Definition  \ref{remark1}.  As $X=\lim_{n\rightarrow +\infty}\inf_{k\geq n} X(k)$, we immediately obtain the conclusion using Definition  \ref{remark1}. 

\noindent{\bf Proof of Lemma \ref{lemma2.2}.}
Set $Y(n):=\inf_{k\geq n}\mathbb{\hat{E}}_{t}[X(k)]$.  Using the arguments analogous to those performed in Lemma \ref{lemma3}, we get $Y(n) \in {L}_{G}^{1^{*}}(\Omega)$.  According to Lemma \ref{lemma2.1}, we obtain
$\lim_{n \rightarrow \infty}\mathbb{\hat{E}}_{t}[Y(n)]=\mathbb{\hat{E}}_{t}[\lim_{n \rightarrow \infty}Y(n)]$.   Because of $Y(n)\leq X(n)$, we derive 
$\mathbb{\hat{E}}_{t}[Y(n)]\leq\mathbb{\hat{E}}_{t}[X(n)]$ and
$
\lim_{n \rightarrow \infty}\mathbb{\hat{E}}_{t}[Y(n)]\leq \varliminf_{n\rightarrow \infty}\mathbb{\hat{E}}_{t}[X(n)],
$
which implies 
$\mathbb{\hat{E}}_{t}[\varliminf_{n\rightarrow \infty}X(n)]\leq \varliminf_{n\rightarrow \infty}\mathbb{\hat{E}}_{t}[X(n)]
$ we expect.

Now, we are in a position to prove the $G$-martingale convergence theorem step-by-step using the uppercrossing inequality.

\begin{definition}\label{upcrossingd}
	A random time $\tau: \Omega\rightarrow [0,+\infty)$ is called an $\mathbb{F}$-stopping time, if $\{\tau \leq t\}\in \mathscr{F}_{t}$ for every $t\geq 0$. 
\end{definition}
\begin{definition}\label{stoppingtime}
	For a finite subset $F\subset [0,+\infty)$, the interval $[\alpha,\beta]$ and the process $M=\{M(t)\}$ with $M(t)\in L_{G}^{1}(\Omega)$, we define the a series of $\mathbb{F}$-stopping times recursively by:
	$$
	\begin{array}{l}
		\tau_{1}(\omega)=\min \left\{t \in F ; M(t;\omega) \leq \alpha\right\},~
		\sigma_{j}(\omega)=\min \left\{t \in F ; t \geq \tau_{j}(\omega), \ \ M(t;\omega)\geq\beta\right\},\\
		\tau_{j+1}(\omega)=\min \left\{t \in F ; t \geq \sigma_{j}(\omega), \ \  M(t;\omega)\leq\alpha\right\}.
	\end{array}
	$$
	And the minimum of an empty set is defined as $+\infty$. Let $U_{F}(\alpha, \beta ; M(\omega))$ be the largest number $j$ such that $\sigma_{j}(\omega)<+\infty$. For any general set $I\subset [0,+\infty)$, we define
	$
	U_{I}(\alpha, \beta ; M(\omega))=\sup \left\{U_{F}(\alpha, \beta ; M(\omega)) ; F \subseteq I, \  F \ {\text is} \text { finite}\right\}.
	$
\end{definition}
\begin{proposition}[Upcrossing Inequality, A Discrete Version, \cite{Li-31}]\label{upcrossing}
	Assume that $\{-M(n):n=1,2,\cdots,N\}$ is a $G$-supermartingale.  If  $M(n)\in L_{G}^{1}(\Omega_{n})$, then we have
	$
	\mathbb{\hat{E}}[U_{\{1,2,\cdots,N\}}(\alpha, \beta ; M(\omega))]\leq \frac{\mathbb{\hat{E}}[(M(N)-\alpha)^{+}]}{\beta-\alpha}.
	$
\end{proposition}

\begin{lemma}[Uppercrossing Inequality, A Continuous Version]\label{upcrossing2}
	Assume that $\{M(t):t\in[0,+\infty)\}$ is a right- or left-continuous function and $\{-M(t):t\in[0,+\infty)\}$ is a $G$-supermartingale.  If $M(t)\in L_{G}^{1}(\Omega_{t})$, then we have that, for any integer $n>0$, 
	$
	\mathbb{\hat{E}}[U_{[0,n]}(\alpha, \beta ; M(\omega))]\leq \frac{\mathbb{\hat{E}}[(M({n})-\alpha)^{+}]}{\beta-\alpha}.
	$
\end{lemma}
\noindent{\bf Proof.}
Define $A_{j}:=\cup_{1\leq k\leq j}\{ni/k: i=0,1,\cdots, k\}$. Then, the monotone convergence theorem (Theorem \ref{monotone}), together with  Definition \ref{upcrossingd} and Proposition \ref{upcrossing},  immediately yields:
$
\mathbb{\hat{E}}
\left[U_{[0,n]\cap \mathbb{Q}}(\alpha, \beta ; M(\omega))\right]=\lim_{j\rightarrow +\infty}	\mathbb{\hat{E}}[U_{A_{j}}(\alpha, \beta ; M(\omega))]\leq \frac{\mathbb{\hat{E}}[(M({n})-\alpha)^{+}]}{\beta-\alpha}.
$
Thus, for any sufficiently small $\epsilon>0$, as $M$ is right- or left-continuous, 
$
\mathbb{\hat{E}}\left[U_{[0,n]}(\alpha, \beta ; M(\omega))\right]\leq 	\mathbb{\hat{E}}\left[U_{[0,n]\cap \mathbb{Q}}(\alpha+\epsilon, \beta-\epsilon ; M(\omega))\right]\leq \frac{\mathbb{\hat{E}}[(M({n})-\alpha)^{+}]}{\beta-\alpha-2\epsilon},
$
which validates the conclusion as required due to the arbitrariness of $\epsilon$'s selection. 

\noindent{\bf Proof of Proposition \ref{convergence}.}
From Lemma \ref{upcrossing2} and Proposition \ref{monotone}, it follows that
\begin{eqnarray*}
	&  & \mathbb{\hat{E}}[U_{[0,+\infty)}(\alpha, \beta ; -M(\omega))]=\lim_{n\rightarrow +\infty}\mathbb{\hat{E}}[U_{[0,n]}(\alpha, \beta ; -M(\omega))]
	\leq \sup_{n\in \mathbb{N}}\frac{\mathbb{\hat{E}}[(-M({n})-\alpha)^{+}]}{\beta-\alpha} \leq \\
	& & \frac{\sup_{t\geq 0}\mathbb{\hat{E}}[(-M)^{+}({t})]+\vert \alpha\vert }{\beta-\alpha}
	=\frac{\sup_{t\geq 0}\mathbb{\hat{E}}[M^{-}({t})]+\vert \alpha\vert }{\beta-\alpha}
	\leq \frac{\mathbb{\hat{E}}[\sup_{t\geq 0}M^{-}({t})]+\vert \alpha\vert }{\beta-\alpha} <+\infty.
\end{eqnarray*}
So $U_{[0,+\infty)}(\alpha, \beta ; -M(\omega))<+\infty$ quasi-surely. Denote by
$
A_{\alpha,\beta}:=\big\{U_{[0,+\infty)}(\alpha, \beta ; -M(\omega))=+\infty\big\}.
$
Since $\{\omega:-M(t;\omega) \ {\rm does} \ {\rm not} \ {\rm converge}\}\subset \cup_{\alpha,\beta\in \mathbb{Q}}A_{\alpha,\beta}$, $-M(t)$ converges quasi-surely to some $-M(+{\infty})$. Here, $M(+{\infty})$ can be $+\infty$ or $-\infty$. By the fact that 
$
M(t)\geq \inf_{t\geq 0}-M^{-}(t)=-\sup_{t\geq 0}M^{-}(t)
$
and Lemma \ref{lemma3}, we have $M(+{\infty})\in \mathcal{L}_{G}^{{1}^{*}_{*}}(\Omega)$. And by Lemma \ref{lemma2.2}, we further have
\begin{eqnarray*}
	\mathbb{\hat{E}}[\vert M(+{\infty})\vert ]&&\leq \varliminf_{n \rightarrow \infty} {\mathbb{\hat{E}}}\left[\left\vert M({n})\right\vert \right]< 2\mathbb{\hat{E}}\left[\sup_{n\in \mathbb{N}}M^{-}({n})\right]+\varliminf_{n \rightarrow \infty}\mathbb{\hat{E}}[M({n})]\\
	&&\leq 2\mathbb{\hat{E}}\left[\sup_{t\geq 0}M^{-}(t)\right]+\mathbb{\hat{E}}[M({1})]<\infty.
\end{eqnarray*}
Thus, $M({+\infty})$, finite quasi-surely, belongs to ${L}_{G}^{{1}^{*}_{*}}(\Omega)$.  Finally,  by virtue of Lemma \ref{lemma2.2},
we have
$
\mathbb{\hat{E}}_{t}[M(+{\infty})]\leq \varliminf _{k\rightarrow +\infty}\mathbb{\hat{E}}_{t}[M({t_{k}})]\leq M({t}),
$
which completes the proof.	
\subsection{Proof of Proposition \ref{exist2}}\label{appendix3}
From Propositions \ref{GIto} and \ref{exist},  it follows that
\begin{eqnarray*}
	&&V({\bm x}(t\wedge \tau_{N}),t\wedge \tau_{N})=V({\bm x}_{0},0)+\int_{0}^{t\wedge \tau_{N}} V_{t}({\bm x}(s),s){\rm d}s\\
	&&+\int_{0}^{t\wedge \tau_{N}} V_{x_{i}}(x(s),s)f^{i}({\bm x}(s),s){\rm d}s  
	+ \int_{0}^{t\wedge \tau_{N}} V_{x_{i}}({\bm x}(s),s)g^{ij}({\bm x}(s),s){\rm d}B_{j}(s)\\
	&&+ \int_{0}^{t\wedge \tau_{N}} V_{x_{k}}({\bm x}(s),s)h^{kij}({\bm x}(s),s){\rm d}\langle B_{i}, B_{j} \rangle(s)\\
	&&+\int_{0}^{t\wedge \tau_{N}}\frac{1}{2} V_{x_{k}x_{l}}({\bm x}(s),s)g^{ki}({\bm x}(s),s)
	g^{lj}({\bm x}(s),s){\rm d}\langle B_{i}, B_{j} \rangle(s).
\end{eqnarray*}
Set ${\bm \eta}=(\kappa_{ij})\in M_{G}^{1}([0,T]; \mathbb{S}^{m})$, where ${ \eta}_{ij}=V_{x_{k}}(h^{kij}+h^{kji})+V_{x_{k}x_{l}}g^{ki}g^{lj}$. Using Proposition \ref{prop5} leads us to the calculations as follows:
\begin{eqnarray*}
	&&V({\bm x}(t\wedge \tau_{N}),t\wedge \tau_{N})=V({\bm x}_{0},0)+\int_{0}^{t\wedge \tau_{N}} V_{t}({\bm x}(s),s){\rm d}s\\
	&&+\int_{0}^{t\wedge \tau_{N}} V_{x_{i}}({\bm x}(s),s)f^{i}({\bm x}(s),s){\rm d}s  
	+ \int_{0}^{t\wedge \tau_{N}} V_{x_{i}}({\bm x}(s),s)g^{ij}({\bm x}(s),s){\rm d}B_{j}(s)\\
	&&+\int_{0}^{t\wedge \tau_{N}} \frac{1}{2}\kappa_{ij}({\bm x}(s),s) {\rm d}\langle B_{i},B_{j} \rangle \\
	&\leq&V({\bm x}_{0},0)+\int_{0}^{t\wedge \tau_{N}} V_{t}({\bm x}(s),s){\rm d}s+\int_{0}^{t\wedge \tau_{N}} G(\bm \eta) {\rm d}s\\
	&&+\int_{0}^{t\wedge \tau_{N}} V_{x_{i}}({\bm x}(s),s)f^{i}({\bm x}(s),s){\rm d}s  
	+ \int_{0}^{t\wedge \tau_{N}} V_{x_{i}}({\bm x}(s),s)g^{ij}({\bm x}(s),s){\rm d}B_{j}(s)
\\
	&=&V({\bm x}_{0},0)+\int_{0}^{t\wedge \tau_{N}} \mathcal{L} V({\bm x}(s),s){\rm d}s 
	+ \int_{0}^{t\wedge \tau_{N}} V_{x_{i}}({\bm x}(s),s)g^{ij}({\bm x}(s),s){\rm d}B_{j}(s)\\
	&\leq &V({\bm x}_{0},0)+\int_{0}^{+\infty} \gamma(t){\rm d}t
	+ \int_{0}^{t\wedge \tau_{N}} V_{x_{i}}({\bm x}(s),s)g^{ij}({\bm x}(s),s){\rm d}B_{j}(s)
\end{eqnarray*}
Then,
$
	\mathbb{\hat{E}}[\vert V({\bm x}(t\wedge \tau_{N}),t\wedge \tau_{N})\vert ]\leq \vert V({\bm x}_{0},0)\vert + \int_{0}^{+\infty} \gamma(t){\rm d}t:=K<+\infty, 
$
which implies
\begin{eqnarray}\label{inequ}
	\infty&>&K\geq \mathbb{\hat{E}}[\vert V({\bm x}(t\wedge \tau_{N}),t\wedge \tau_{N})\vert ]\geq \mathbb{\hat{E}}[\mu(\vert {\bm x}(t\wedge \tau_{N})\vert)]\geq   \nonumber\\
	&\geq& \mu(N)c(\tau_{N}\leq t)\geq  \mu(N)c(\tau_{\infty}\leq t)
\end{eqnarray}
where $\mu(r):=\inf_{\vert \bm x\vert \geq r, t\geq 0}V(\bm x,t)$ and $\lim_{r\rightarrow+\infty}\mu(r)=+\infty$ because of the  condition assumed in \eqref{conditionofV2}.
Now, letting $N\rightarrow+\infty $ in \eqref{inequ} yields $c(\tau_{\infty}\leq t)=0$ for any $t$.  Finally, further letting $t\rightarrow+\infty $ gives $c(\tau_{\infty}\leq +\infty)=0$, which completes the proof.

%Here for each time step and each $i=1,2, \cdots, m$, we select a random $\sigma_{i}(t_{n})$ from the interval $[\sigma_{r}, \sigma_{s}]$ uniformly and independently. In addition, $\Delta{B}(t_{n}) \sim \mathcal{N}(0,\sigma_{i}(t_{n}) \Delta t)$, $\Delta \langle B_{i} \rangle (t_n)=(\sigma_{i}(t_{n}))^{2}{\Delta t}$. That is, the probability itself has an uncertainty for each time step. For all these probability measures, they satisfy the same conclusion if it is correct quasi-surely.

\subsection{Proof of Lemma \ref{theorembarbalat}}\label{appendix_barbalat}

\textcolor{black}{To prove Lemma \ref{theorembarbalat}, we first establish the inequality as follows.}

\begin{lemma} \label{theoreminequality}
	For $A_{ij}(t)\in M_{G}^{2}[0, T]$, denote by ${\bm A}(t)=(a_{ij}(t))_{d\times m}$.  Then, we have
	$
	\mathbb{\hat{E}}\left\vert \int_{0}^{T} {\bm A}(t) \mathrm{~d} {\bm B}(t)\right \vert ^2 \leq d\bar{\gamma} \  \mathbb{\hat{E}}\int_{0}^{T} \vert {\bm A}(t)\vert ^{2} \mathrm{~d} t.
	$
\end{lemma}

\noindent{\bf Proof.} For simplicity of expression, we apply Einstein's notations \cite{einstein1922general} in the following arguments and throughout if they are necessary.   From Theorem \ref{BDG} and Remark \ref{lemmaqua}, it follows that
\begin{eqnarray*}
	&&\mathbb{\hat{E}}\left\vert \int_{0}^{T} {\bm A}(t) \mathrm{~d} {\bm B}(t) \right\vert ^2 =\mathbb{\hat{E}}  \left (\int_{0}^{T}a_{ij}(t) \mathrm{~d} B_{j}(t) \int_{0}^{T}a_{ik}(t) \mathrm{~d} B_{k}(t) \right )\\
	&=& \mathbb{\hat{E}} \int_{0}^{T}a_{ij}(t)a_{ik}(t) \mathrm{~d} \langle B_{j}, B_{k}\rangle(t)
	=\mathbb{\hat{E}} \int_{0}^{T} {\rm trace}({\bm A}(t)\mathrm{~d} \langle {\bm B} \rangle(t) {\bm A}^{\top}(t) )\\
	&\leq & d\cdot \mathbb{\hat{E}} \int_{0}^{T} \lambda_{\rm max}(A(t)\mathrm{~d} \langle B \rangle(t) A^{\top}(t) )
	= d\cdot\mathbb{\hat{E}} \int_{0}^{T} \vert {\bm A}(t)\mathrm{~d} \langle {\bm B} \rangle(t) {\bm A}^{\top}(t) \vert_{2} \\
	&=& d\cdot\mathbb{\hat{E}} \int_{0}^{T} \vert {\bm A}(t)\vert_{2} ^{2} \mathrm{~d} \vert \langle {\bm B} \rangle\vert_{2} (t)
	\leq  d\bar{\gamma} \cdot \mathbb{\hat{E}} \int_{0}^{T} \vert {\bm A}(t)\vert_{2} ^{2} \mathrm{~d} t
	\leq  d\bar{\gamma} \cdot \mathbb{\hat{E}} \int_{0}^{T} \vert {\bm A}(t)\vert^{2} \mathrm{~d} t.
\end{eqnarray*}
The proof is therefore completed.	

\noindent{\bf Proof of Lemma \ref{theorembarbalat}.} Now, we need to prove the lemma using contradiction.  If this is not true, then there exists $Q \in \mathcal{Q} $ such that
$Q \left(\left\{\omega :\liminf_{t \rightarrow +\infty} \eta(\bm x(t;\omega))>0\right\} \right)>0.
$
Thus, there exists $\epsilon>0$ such that $Q(\Omega_{1})\geq 2\epsilon$ with
$
\Omega_{1}=\left\{\omega\in \Omega_{0}: \liminf_{t \rightarrow +\infty} \eta(\bm x(t))>2\epsilon\right\}.
$
Since
$
\Omega_{1}=\cup_{n=1}^{+\infty} \left( \Omega_{1} \cap \left\{\omega: \sup_{t \geq 0}\vert \bm x(t;\omega)\vert <n \right\}\right),
$
there exists a number $N>0$ such that $Q(\Omega_{2})\geq \epsilon$ in which
$
\Omega_{2}= \Omega_{1} \cap \left\{\omega: \sup_{t \geq 0}\vert \bm x(t;\omega)\vert <N \right\}$.

Now, we define the $\mathbb{F}$-stopping times as 
$$
\begin{array}{l}
	\sigma_{1}(\omega):=\inf\{t: \eta(\bm x(t;\omega))\geq 2\epsilon\}, ~~
	\sigma_{2i}(\omega):=\inf\{t: \eta(\bm x(t;\omega))\leq \epsilon, ~ t\geq \sigma_{2i-1}(\omega)\}, \\
	\sigma_{2i+1}(\omega):=\inf\{t: \eta(\bm x(t;\omega))\geq 2\epsilon, ~ t\geq \sigma_{2i}(\omega)\}, ~~
	\tau_{N}(\omega):=\inf\{t: \vert \bm x(t;\omega)\vert\geq N\}.
\end{array}
$$
For all $\omega \in \Omega_{2}$, $\tau_{N}(\omega)=+\infty$ and $\sigma_{i}(\omega)<+\infty$ for all $i>0$ using the formula \eqref{integralbound} and the definition of $\Omega_1$. By virtue of Proposition \ref{exist}, ${\bm M}(t)=\int_{0}^{t\wedge \tau_{N}}{\bm g}({\bm x}(s),s) {\rm d}{\bm B}(s)$ is a $Q$-martingale for each $Q \in \mathcal{Q}$.   Hence, using Assumption \ref{asump1}, Lemma \ref{theoreminequality}, Hölder's inequality, and Doob's martingale inequality in traditional stochastic analysis, we obtain that for each $T>0$,
\begin{eqnarray*}
	&&	E_{Q}[1_{\{\tau_{N}\wedge \sigma_{2i-1}<+\infty\}}\sup_{0 \leq t \leq T}\vert \bm x(\tau_{N}\wedge (\sigma_{2i-1}+t))-\bm x(\tau_{N}\wedge \sigma_{2i-1})\vert ^{2}]
	\\
	&\leq& 3E_{Q}\left [1_{\{\tau_{N}\wedge \sigma_{2i-1}<+\infty\}}\sup_{0 \leq t \leq T}\left \vert \int_{\tau_{N}\wedge \sigma_{2i-1}}^{\tau_{N}\wedge (\sigma_{2i-1}+t)}  \bm f(\bm x(s),s) {\rm d}s   \right \vert ^{2}\right ] \\
	& &+3E_{Q}\left [1_{\{\tau_{N}\wedge \sigma_{2i-1}<+\infty\}}\sup_{0 \leq t \leq T}\left \vert \int_{\tau_{N}\wedge \sigma_{2i-1}}^{\tau_{N}\wedge (\sigma_{2i-1}+t)}  \bm g(\bm x(s),s) {\rm d}\bm B(s)   \right \vert ^{2}\right ] \\
	& &+3E_{Q}\left [1_{\{\tau_{N}\wedge \sigma_{2i-1}<+\infty\}}\sup_{0 \leq t \leq T}\left \vert \int_{\tau_{N}\wedge \sigma_{2i-1}}^{\tau_{N}\wedge (\sigma_{2i-1}+t)}  \bm h(\bm x(s),s) {\rm d}\langle \bm B \rangle(s)   \right \vert ^{2}\right ] \\
	&\leq& 3TE_{Q}\left [1_{\{\tau_{N}\wedge \sigma_{2i-1}<+\infty\}}\sup_{0 \leq t \leq T}\int_{\tau_{N}\wedge \sigma_{2i-1}}^{\tau_{N}\wedge (\sigma_{2i-1}+t)}  \vert \bm f(\bm x(s),s)\vert ^{2} {\rm d}s  \right ] \\
	& &+12E_{Q}\left [1_{\{\tau_{N}\wedge \sigma_{2i-1}<+\infty\}}\left \vert \int_{\tau_{N}\wedge \sigma_{2i-1}}^{\tau_{N}\wedge (\sigma_{2i-1}+T)}   \bm g(\bm x(s),s)  {\rm d} \bm B(s) \right \vert \right ] \\
	& &+3T\bar{\gamma}^{2}m^2E_{Q}\left [1_{\{\tau_{N}\wedge \sigma_{2i-1}<+\infty\}}\sup_{0 \leq t \leq T}\int_{\tau_{N}\wedge \sigma_{2i-1}}^{\tau_{N}\wedge (\sigma_{2i-1}+t)}  \vert \bm h(\bm x(s),s)\vert ^{2} {\rm d}s  \right ]
\end{eqnarray*}
\begin{eqnarray*}	
	&\leq& 3T\mathbb{\hat{E}}\left [1_{\{\tau_{N}\wedge \sigma_{2i-1}<+\infty\}}\int_{\tau_{N}\wedge \sigma_{2i-1}}^{\tau_{N}\wedge (\sigma_{2i-1}+T)}  \vert \bm f(\bm x(s),s)\vert ^{2} {\rm d}s  \right ]\\
	& &+12 d \bar{\gamma}\mathbb{\hat{E}}\left [1_{\{\tau_{N}\wedge \sigma_{2i-1}<+\infty\}}\int_{\tau_{N}\wedge \sigma_{2i-1}}^{\tau_{N}\wedge (\sigma_{2i-1}+T)}  \vert \bm g(\bm x(s),s)\vert ^{2} {\rm d}s  \right ] \\
	& &+3T\bar{\gamma}^{2}m^2\mathbb{\hat{E}}\left [1_{\{\tau_{N}\wedge \sigma_{2i-1}<+\infty\}}\int_{\tau_{N}\wedge \sigma_{2i-1}}^{\tau_{N}\wedge (\sigma_{2i-1}+T)}  \vert \bm h(\bm x(s),s)\vert ^{2} {\rm d}s  \right ]\\
	&\leq& 3K_{N}^{2} T(T+4d\bar{\gamma}+T\bar{\gamma}^{2}m^2).
\end{eqnarray*}
As $\eta$ is continuous, there exists a number $\delta>0$ such that, for every $\bm x,\bm y \in B(N)$ and $\vert \bm x-\bm y\vert \leq \delta$, $\vert \eta(\bm x)-\eta(\bm y)\vert <\epsilon$.
We select sufficiently small $T>0$ such that
$
{ 3K_{N}^{2}T(T+4d\bar{\gamma}+T\bar{\gamma}^{2}m^2)}/{\delta^2}<\frac{\epsilon}{2}.
$
Thus, we have
$Q\left(1_{\{\tau_{N}\wedge \sigma_{2i-1}<+\infty\}}\sup_{0 \leq t \leq T}\vert \bm x(\tau_{N}\wedge (\sigma_{2i-1}+t))-\bm x(\tau_{N}\wedge \sigma_{2i-1})\vert  \geq \delta \right)\leq { 3K_{N}^{2}T(T+4d\bar{\gamma}+T\bar{\gamma}^{2}m^2)}/{\delta^2}<\frac{\epsilon}{2}$.
Hence, we have
$Q(\{\sigma_{2i-1}<+\infty, \tau_{N}=+\infty\}\cap \{\sup_{0\leq t\leq T}\vert \bm x(\sigma_{2i-1}+t)-\bm x(\sigma_{2i-1})\vert \geq \delta\}) \leq \frac{\epsilon}{2}.
$
By the definition and the property 
of $\Omega_{2}$, we conclude that
$
Q\left(\{\sigma_{2i-1}<+\infty, \tau_{N}=+\infty\}\cap\left\{\sup_{0\leq t\leq T}\vert \bm x(\sigma_{2i-1}+t)-\bm x(\sigma_{2i-1})\vert < \delta\right\}\right)
\geq \epsilon-\frac{\epsilon}{2}=\frac{\epsilon}{2},
$
which further implies that
\begin{eqnarray*}
	&&Q\left(\{\sigma_{2i-1}<+\infty, \tau_{N}=+\infty\}\cap\left\{\sup_{0\leq t\leq T}\vert \eta(\bm x(\sigma_{2i-1}+t))-\eta(\bm x(\sigma_{2i-1}))\vert < \epsilon\right\}\right)\\
	&&
	\geq Q\left(\{\sigma_{2i-1}<+\infty, \tau_{N}=+\infty\}
	\cap\left\{\sup_{0\leq t\leq T}\vert \bm x(\sigma_{2i-1}+t)-\bm x(\sigma_{2i-1})\vert < \delta\right\}\right)
	\geq \frac{\epsilon}{2}.	
\end{eqnarray*}
Define 
$
\tilde{\Omega}_{i}:=\left  \{ \sup_{0\leq t \leq T} \vert \eta(\bm x(\sigma_{2i-1}+t))-\eta(\bm x (\sigma_{2i-1}))\vert < \epsilon \right \}.
$
Then, on $\tilde{\Omega}_{i}\cap \{\sigma_{2i-1}<+\infty\}$, we have $\sigma_{2i}-\sigma_{2i-1}\geq T$.   By \eqref{limitofeta}, if $\sigma_{2i-1}<+\infty$, then $\sigma_{2i}<+\infty$ quasi-surely. Thus, 
\begin{eqnarray*}
	+\infty &>& \mathbb{\hat{E}} \int_{0}^{+\infty} \eta(\bm x(t)){\rm d}t \geq	E_{Q} \int_{0}^{+\infty} \eta(\bm x(t)){\rm d}t\\
	&\geq &\sum_{i=1}^{+\infty} E_{Q}\left [1_{\{\tau_{N}=+\infty, \sigma_{2i-1}<+\infty,\sigma_{2i}<+\infty\}} \int_{\sigma_{2i-1}}^{\sigma_{2i}} \eta(\bm x(t)){\rm d}t\right ]\\
	&\geq& \epsilon \sum_{i=1}^{+\infty} E_{Q}\left [1_{\{\tau_{N}=+\infty, \sigma_{2i-1}<+\infty\}} (\sigma_{2i}-\sigma_{2i-1})\right ]\\
	&\geq& \epsilon \sum_{i=1}^{+\infty} E_{Q}\left [1_{\{\tau_{N}=+\infty, \sigma_{2i-1}<+\infty\}\cap \tilde{\Omega}_{i}} (\sigma_{2i}-\sigma_{2i-1})\right ]\\
	&\geq& \epsilon T  \sum_{i=1}^{+\infty} Q(\{\tau_{N}=+\infty, \sigma_{2i-1}<+\infty\}\cap \tilde{\Omega}_{i})
	\geq \epsilon T  \sum_{i=1}^{+\infty} \frac{\epsilon}{2} =+\infty,
\end{eqnarray*}
which indicates a contradiction.  Consequently, we get $\lim_{t \rightarrow +\infty}\eta(\bm x(t))=0$ quasi-surely. 

\subsection{Dynamic Stability in Example \ref{exstability}}\label{app4}
Here, we validate the quasi-sure stability of the considered equations in Example \ref{exstability}.   To this end, we set $V(\bm x):=\vert \bm x\vert^{\alpha}$ for some given $0<\alpha<1$, which yields
$
	\mathcal{L}V(\bm x)= \alpha \vert \bm x \vert^{\alpha-2} \langle \bm x, \bm f (\bm x)\rangle+G\left(\left(k^2(\alpha-1)\alpha \vert \bm x \vert^{\alpha}\right)_{i,j=1}^{m}\right),
$
where $\left(k^2(\alpha-1)\alpha \vert \bm x \vert^{\alpha}\right)_{i,j=1}^{m}$ stands for an $m\times m$ matrix such that all elements are $k^2(\alpha-1)\alpha \vert \bm x \vert^{\alpha}$. As $c_{-1}:=(-1)_{i,j=1}^{m}$ is a non-positive symmetric matrix with eigenvalues 0 and $-m$, we have $c_{-1}<0$. Set $0<\alpha<1+\frac{L}{k^2c_{-1}}<1$, we obtain that
$
	\mathcal{L}V(\bm x)= \alpha \vert \bm x \vert^{\alpha-2} \langle \bm x, \bm f (\bm x)\rangle+k^2c_{-1}(1-\alpha)\alpha \vert \bm x \vert^{\alpha}
	\leq\alpha \vert \bm x \vert^{\alpha}(L+k^2c_{-1}(1-\alpha)).
$
Set $\eta({\bm x}):=\alpha \vert \bm x \vert^{\alpha}(L+k^2c_{-1}(1-\alpha))<0$.  Hence,  in light of Proposition \ref{exist2}  and Theorem \ref{lasalle1}, if we could confirm a \textit{statement} that the system in Example $\ref{exstability}$ does not reach $\bm 0$ before it explodes, $V(\bm x)$ with $\alpha<1$ and along any trajectory apart from $\bm 0$ is differentiable to the second order, so that the quasi-sure convergence of $\bm x$ is guaranteed to $\bm 0$, the kernel of $\eta$.   To make confirm the statement, we first introduce the following result.

\begin{proposition}\label{prop5ex}
	Let 
	$
	M(t)=\int_{0}^{t} \kappa_{ij}(s){\rm d}\langle B_{i},B_{j} \rangle (s)+\int_{0}^{t} 2G(-\bm \eta){\rm d}s,
	$
	where ${\bm \eta} \in M_{G}^{1}([0,T];\mathbb{S}^{m})$. Then, we have $M(t)\geq 0$ quasi-surely. Particularly, $\mathbb{\hat{E}}[M(t)]\geq 0$. 
\end{proposition}

The proof of the above proposition is akin to the proof for Proposition \ref{prop5}, which is omitted here.

Now, we make the final confirmation.  We set $\tau_{N}:=\inf \{t\geq 0:\vert {\bm x(t)}\vert \geq N\}$ and $\xi_{\epsilon}=\inf \{t\geq 0:\vert {\bm x(t)}\vert \leq \epsilon\}$ for $\epsilon,N>0$, and select $V(\bm x)=\log\vert \bm x \vert$. Then,  using the formula presented in Theorem~\ref{GIto} and Proposition~\ref{exist}, we get
\begin{eqnarray*}
	\log \vert \bm x(t\wedge \tau_{N} \wedge \xi_{\epsilon}) \vert&=&\log \vert \bm x_{0} \vert +\int_{0}^{t\wedge \tau_{N} \wedge \xi_{\epsilon}}\frac{\langle \bm x(s), \bm f(\bm x(s)) \rangle}{\vert \bm x \vert^2}\rm d s\\
	&&+\sum_{j=1}^{n}\int_{0}^{t\wedge \tau_{N} \wedge \xi_{\epsilon}}k{\rm d} B_{j}(s)-\sum_{i,j=1}^{n}\int_{0}^{t\wedge \tau_{N} \wedge \xi_{\epsilon}}\frac{1}{2}k^2 \langle B_{i},B_{j} \rangle(s)
\end{eqnarray*}
Noticing the local Lipschitz property of $\bm f$ gives $\vert \langle \bm x, \bm f(\bm x) \rangle\vert \leq \vert \bm x \vert \vert \bm f(\bm x) \vert\leq K_{N}\vert \bm x \vert^{2} $ on $[0,\tau_{N})$. Set $c_{1}:=G((1)_{i,j=1}^{m})>0$. Then, by Proposition \ref{prop5ex},  we have
$
	\mathbb{\hat{E}}[\log \vert \bm x(t\wedge \tau_{N} \wedge \xi_{\epsilon}) \vert\geq \mathbb{\hat{E}}[\log \vert \bm x_{0} \vert]-\int_{0}^{t\wedge \tau_{N} \wedge \xi_{\epsilon}}(K_{N}+k^2c_{1}){\rm d}s]
	\geq  \mathbb{\hat{E}}[\log \vert \bm x_{0}\vert]-(K_{N}+k^2c_{1})t.
$
On the other hand,
$
	\mathbb{\hat{E}}[\log \vert \bm x(t\wedge \tau_{N} \wedge \xi_{\epsilon}) \vert ]
	\leq c(\xi_{\epsilon}< t\wedge \tau_{N} )\log\epsilon + c(\xi_{\epsilon}\geq t\wedge \tau_{N} )\log N
	\leq  c(\xi_{\epsilon}< t\wedge \tau_{N} )\log\epsilon + \log N.
$
Hence, we obtain 
$
	\mathbb{\hat{E}}[\log \vert \bm x_{0}\vert]-(K_{N}+k^2c_{1})t \leq c(\xi_{\epsilon}< t\wedge \tau_{N} )\log\epsilon + \log N.
$
First, letting $\epsilon \rightarrow 0$ results in $c(\xi_{0}< t\wedge \tau_{N} )=0$.  Then, letting both $t$ and $N \rightarrow +\infty$ yields $c(\xi_{0}< \tau_{\infty} )=0$, which confirms the above statement and finally completes the proof. 

\subsection{Invariant Set Associated with Autonomous $G$-SDEs}\label{appendix2}

\begin{theorem}\label{app1}
	We consider the following autonomous $G$-SDEs:
	\begin{equation}\label{equationapp1}
		{\rm d}{\bm x}(t)={\bm f}({\bm x}(t)){\rm d}t+{\bm g}({\bm x}(t)){\rm d}{\bm B}(t)+{\bm h}({\bm x}(t)){\rm d}\langle {\bm B} \rangle (t),
	\end{equation}
	where ${\bm f}:\mathbb{R}^{d}\rightarrow \mathbb{R}^{d}$, ${\bm g}:\mathbb{R}^{d} \rightarrow \mathbb{R}^{d\times m}$, ${\bm h}:\mathbb{R}^{d} \rightarrow \mathbb{R}^{d\times m^2}$, and ${\bm f}({\bm a})={\bm g}({\bm a})={\bm h}({\bm a})={\bm 0}$.  Clearly, ${\bm f},{\bm g}$ and ${\bm h}$ are all globally Lipschitzian. Then, we have that, for all ${\bm x}_{0}\neq a$,
$
		c\big(\{\omega: \exists~t>0, ~~\bm{x}(t,\omega;{\bm x}_{0})={\bm a}\} \big)=0,
$
	which indicates that the trajectory does not approach ${\bm a}$ quasi-surely in a finite time.	
\end{theorem}

\noindent{\bf Proof.}
We know that the $G$-SDEs \eqref{equationapp1} have a unique solution on $M_{G}[0,T] $ for every $T>0$ according to \cite{Peng-13}. 
First, we need to perform the proof for the situation of $\bm a={\bm 0}$.
Now set 
$
\mathcal{A}:=\{\omega: {\bm x}(t,\omega)={\bm 0} \  {\rm for} \ {\rm some} \ t \in [0,+\infty)\}.
$
If $c(\mathcal{A})>0$, then there exists a number $T>0$ such that $c(\mathcal{A}_{T})>0$ where
$
	\mathcal{A}_{T}:=\{\omega: {\bm x}(t,\omega)={\bm 0} \  {\rm for} \ {\rm some} \ t \in [0,T]\},
$
which is due to the fact that $\mathcal{A}=\cup_{T=1}^{+\infty}\mathcal{A}_{T}$.  Next, introduce the stopping time 
$
	\tau_{\epsilon}:=\inf\{t \in [0,+\infty): \vert {\bm x}(t,\omega)\vert \leq \epsilon  \} . 
$
Set $V({\bm x}):=1/\vert {\bm x}\vert =(\vert {\bm x}\vert ^{2})^{-\frac{1}{2}}$.  Then, we perform the calculations using $G$-Itô's formula, obtaining that
\begin{eqnarray*}
&&	V({\bm x}(T\wedge \tau_{\epsilon}))=V({\bm x}_{0})+\int_{0}^{T\wedge \tau_{\epsilon}} V_{x_{i}}({\bm x}(s))f^{i}({\bm x}(s)){\rm d}s  \\
	&&+ \int_{0}^{T\wedge \tau_{\epsilon}} V_{x_{i}}({\bm x}(s))g^{ij}({\bm x}(s)){\rm d}{ B}_{j}(s)
	+\int_{0}^{T\wedge \tau_{\epsilon}} \frac{1}{2}\kappa_{ij}({\bm x}(s)) {\rm d}\langle B_{i},B_{j} \rangle(s)
\\
	&=&V({\bm x}_{0})-\int_{0}^{T\wedge \tau_{\epsilon}} \frac{\langle {\bm x}(s),f({\bm x}(s))\rangle}{\vert {\bm x}\vert ^{3}}{\rm d}s
	-\int_{0}^{T\wedge \tau_{\epsilon}} \frac{x_{i}(s)g^{ij}({\bm x}(s))}{\vert {\bm x}\vert ^{3}}{\rm d}B_{j}(s)\\
	&&+\int_{0}^{T\wedge \tau_{\epsilon}} \Bigg [-\frac{g^{\mu{i}}({\bm x}(s))g^{\mu{j}}({\bm x}(s))}{2\vert {\bm x}\vert ^{3}}
	+\frac{3}{2\vert {\bm x}\vert ^{5}}x_{\mu}x_{v} g^{\mu{i}}({\bm x}(s))g^{\nu{j}}({\bm x}(s))\\
	&&-\frac{x_{v}h^{vij}({\bm x}(s))}{\vert {\bm x}\vert ^{3}}\Bigg ]{\rm d}\langle B_{i},B_{j} \rangle (s)
	\leq  V({\bm x}_{0})+\int_{0}^{T\wedge \tau_{\epsilon}}\Bigg [\frac{\vert {\bm f}({\bm x})\vert }{\vert {\bm x}\vert ^{2}} +\frac{d\bar{\gamma} \vert {\bm g}({\bm x})\vert ^{2}}{2\vert {\bm x}\vert ^{3}}+\frac{3\bar{\gamma}\vert {\bm g}({\bm x})\vert ^{2}}{2\vert {\bm x}\vert ^3}\\
	&&+\frac{\vert {\bm h}({\bm x})\vert \bar{\gamma}}{\vert {\bm x}\vert ^{2}}\Bigg ]{\rm d}s
	+ \int_{0}^{T\wedge \tau_{\epsilon}} V_{x_{i}}({\bm x}(s))g^{ij}({\bm x}(s)){\rm d}B_{j}(s),
\end{eqnarray*}
where $\kappa_{ij}=V_{x_{k}}(h^{kij}+h^{kji})+V_{x_{k}x_{l}}g^{ki}g^{lj}$ and Einstein's notations are applied here.
Let 
$
	\rho({\bm x}):=\frac{\vert {\bm f}({\bm x})\vert }{\vert {\bm x}\vert } +\frac{(d+3)\bar{\gamma} \vert {\bm g}({\bm x})\vert ^{2}}{2\vert {\bm x}\vert ^{2}}+\frac{\vert {\bm h}({\bm x})\vert \bar{\gamma}}{\vert {\bm x}\vert }.
$
Then, there exists a number $K>0$ such that  $\rho ({\bm x})\leq K<+\infty$ because ${\bm f}$, ${\bm g}$ and ${\bm h}$ are globally  Lipschitzian as mentioned above.  Hence,   it follows that 
\begin{eqnarray*}
&&	V({\bm x}(T\wedge \tau_{\epsilon})) \leq V({\bm x}_{0})+ \int_{0}^{T\wedge \tau_{\epsilon}} V({\bm x}(s)) \rho({\bm x}(s)){\rm d}s+ \int_{0}^{T\wedge \tau_{\epsilon}} V_{x_{i}}({\bm x}(s))g^{ij}({\bm x}(s)){\rm d}B_{j}(s)\\
	& & =V({\bm x}_{0})+ \int_{0}^{T} V({\bm x}(s)) \rho ({\bm x}(s))1_{[0,\tau_{\epsilon}]}{\rm d}s
	+ \int_{0}^{T} V_{x_{i}}({\bm x}(s))g^{ij}({\bm x}(s))1_{[0,\tau_{\epsilon}]}{\rm d}B_{j}(s)\\
	&& \leq  V({\bm x}_{0})+K\int_{0}^{T} V({\bm x}(s)) 1_{[0,\tau_{\epsilon}]} {\rm d}s
	+ \int_{0}^{T} V_{x_{i}}({\bm x}(s))g^{ij}({\bm x}(s)) 1_{[0,\tau_{\epsilon}]}{\rm d}B_{j}(s),
\end{eqnarray*}
which implies that
$
	\mathbb{\hat{E}}[V({\bm x}(T\wedge \tau_{\epsilon}))] \leq \mathbb{\hat{E}}[V({\bm x}_{0})]+ K \mathbb{\hat{E}}\int_{0}^{T} V({\bm x}(s)) 1_{[0,\tau_{\epsilon}]} {\rm d}s\\
	\leq \mathbb{\hat{E}}[V({\bm x}_{0})]+ K \int_{0}^{T}\mathbb{\hat{E}}[ V({\bm x}(s\wedge \tau_{\epsilon}))]{\rm d}s.
$
Now, using Gronwall's inequality, we have
$
	\mathbb{\hat{E}}\left [\frac{1}{\vert {\bm x}(T\wedge \tau_{\epsilon})\vert }\right ]\leq \mathbb{\hat{E}}[V({\bm x}_{0})]{\rm e}^{KT}.
$
From the definition of $\tau_{\epsilon}$ and also from the continuity of ${\bm x}(t)$,  it follows that 
$\vert {\bm x}(T\wedge \tau_{\epsilon})\vert =\epsilon$ on $\mathcal{A}_{T}$. Thus,
$
	c(\mathcal{A}_{T})= \epsilon \mathbb{\hat{E}}\left [\frac{1}{\vert {\bm x}(T\wedge \tau_{\epsilon})\vert }1_{\mathcal{A}_{T}} \right ] 
	\leq \epsilon \mathbb{\hat{E}}[V({\bm x}_{0})]{\rm e}^{KT},
$
which is valid for every $\epsilon>0$.  Therefore, we immediately obtain $c(\mathcal{A}_{T})=0$, which is a contradiction.

For the general situation of $\bm{a}$, we set ${\bm y}(t):={\bm x}(t)-{\bm a}$.  Then, $\bm{y}(t)$ satisfies the $G$-SDEs:
$
	{\rm d}{\bm y}(t)={\bm f}({\bm y}(t)+\bm a){\rm d}t+{\bm g}({\bm y}(t)+{\bm a}){\rm d}{\bm B}(t)+{\bm h}({\bm y}(t)+{\bm a}){\rm d}\langle {\bm B} \rangle (t).
$
Consequently, we know that ${\bm y}(t)$ never approaches ${\bm 0}$ quasi-surely, i.e., ${\bm x}(t)$ never approaches ${\bm a}$ quasi-surely.
Therefore, the proof is complete. rew

\subsection{Numerical evidences}\label{appendix_numerical}

Here, we describe the numerical scheme that we use for partially illustrating the analytical results obtained in the main text.   Actually, we do not provide a complete simulation for the solutions of $G$-SDEs  but only 
simulate the corresponding SDEs under a group of probability measures.  A rigorous and complete scheme for simulating the solution of $G$-SDEs still awaits further investigations.

To this end, we first suppose $\bm W(t)$ to be a standard $m$-dimensional Brownian motion on the probability space 
$(\Omega, \mathcal{B}(\Omega), P)$.   Also suppose that $\Theta$ is a bounded, closed and convex subset of $\mathbb{R}^{m\times m}$, where $\Theta=[\underline{\sigma}, \overline{\sigma}]$ for $m=1$.  In addition, 	
$
\mathcal{\tilde{Q}}:=\Big\{P_{\bm \theta} \in \mathcal{M}:P_{\bm \theta} \ {\rm is} \ {\rm the} \ {\rm law} \ {\rm of} \ {\rm process}
 \ \int_{0}^{t}{\bm \theta}(s){\rm d}{\bm W}(s)\ {\rm for~} \forall  \ t\geq 0, {\bm \theta} \in \mathscr{A}_{0,\infty}^{\Theta}\Big\} \subset \mathcal{Q},
$
where $\mathscr{A}_{0,\infty}^{\Theta}$ denotes the collection of all $\Theta$-valued $\mathscr{F}$ adapted function in $[0,+\infty)$. According to Remark 15 in Ref.~\cite{HuPeng-29}, the capacity satisfies $c(\mathcal{A})=\sup_{Q \in \mathcal{\tilde{Q}}}P[\mathcal{A}]$ for any $\mathcal{A}\in \mathscr{B}(\Omega)$, so we can check whether an event is correct quasi-surely on the probability measures space $\mathcal{\tilde{Q}}$.  Thus, we make our numerical simulations on a finite subset of $\mathcal{\tilde{Q}}$ repeatedly as follows and use the case where $\langle B_{i}, B_{j}\rangle=0$ for each $i \neq j$ and all $B_{i}$ are identically distributed.

For the time interval $[0, T]$, we introduce a uniform time partition $0=t_0<t_1<\cdots<t_N=T$ with $\Delta t:=t_{n+1}-t_n=T/N$. We use the following Euler-Maruyama scheme, as proposed in \cite{maruyama1954transition}, to investigate the solution of the SDEs correspondingly from the $G$-SDEs in \eqref{GSDE}:
\begin{equation}\label{Euler}
	{\bm X}({n+1})={\bm X}(n)+{\bm f}({\bm X}(n),t_{n}){\rm \Delta}t+{\bm g}({\bm X}(n),t_{n}){\Delta\bm B}(t_{n})+{\bm h}({\bm X}(n),t_{n})\Delta\langle {\bm B}\rangle (t_{n})
\end{equation}
with $\bm{X}(0)={\bm x}_0$ and $n=0,1,\cdots,N-1$.  \textcolor{black}{Here, $\Delta{B}_i(t_{n}) \sim \mathcal{N}(0,\sigma_{i,n}^{2} \Delta t)$ and $\Delta \langle B_{i} \rangle (t_n)=\sigma_{i,n}^{2}{\Delta t}$ with $\sigma_{i,n}\in [\underline{\sigma}, \overline{\sigma}]$ and $i=0,1,\cdots,m$.}

\textcolor{black}{In order to investigate the dynamics of the corresponding SDEs on 
the probability measures space 
$\mathcal{\tilde{Q}}$, the covariance $\{\sigma_{i,n}\}_{1\leq i\leq m, 1\leq n\leq N}$ should be taken from all the element of the set $[\underline{\sigma}, \overline{\sigma}]^{  m \times N}$.   To do this numerically, we introduce a uniform interval partition $\underline{\sigma}=\sigma_{0}<\sigma_{1}<\cdots<\sigma_{k}=\overline{\sigma}$ with $\Delta\sigma=\sigma_{i+1}-\sigma_{i}=(\overline{\sigma}-\underline{\sigma})/k$. Denote by  $\Sigma_{jl}:=\{ \sigma_i | j\leq i\leq l  \}$, where $1\leq j\leq l \leq k$. For any given tuple $(j,l)$, we choose an element $(\mu_{in})_{1\leq i\leq m, 1\leq n\leq N}\in\Sigma_{jl}^{m\times N}$, set $\sigma_{i,n}=\mu_{in}$ for all $1\leq i\leq m, 1\leq n\leq N$, and then approximate the dynamics of the SDEs correspondingly from \eqref{GSDE} using the scheme specified in \eqref{Euler}, which enables us to numerically produce a large number of simulating trials. }

\begin{figure}[htbp]
	\centering
 \subfigure{
		\includegraphics[width=4cm]{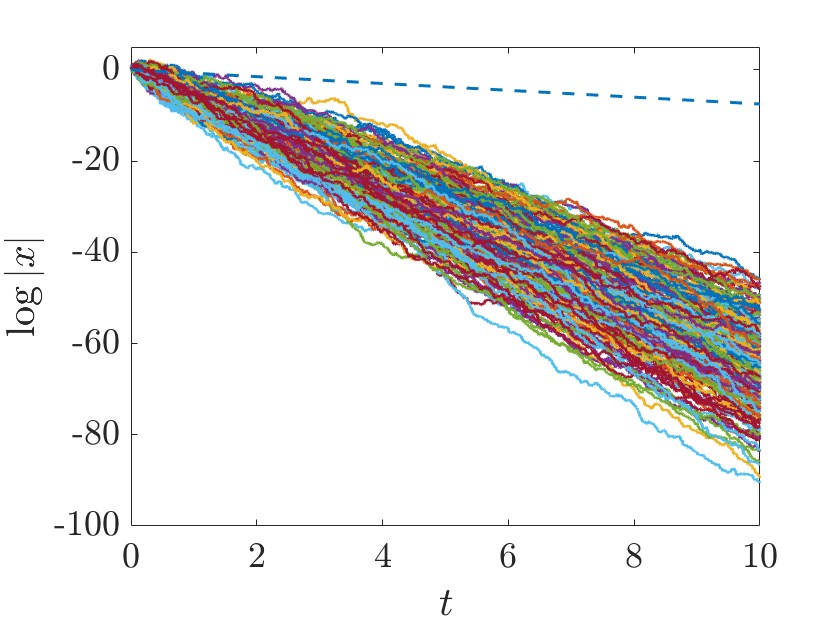}
		\label{1a}}
  \subfigure{
		\includegraphics[width=4cm]{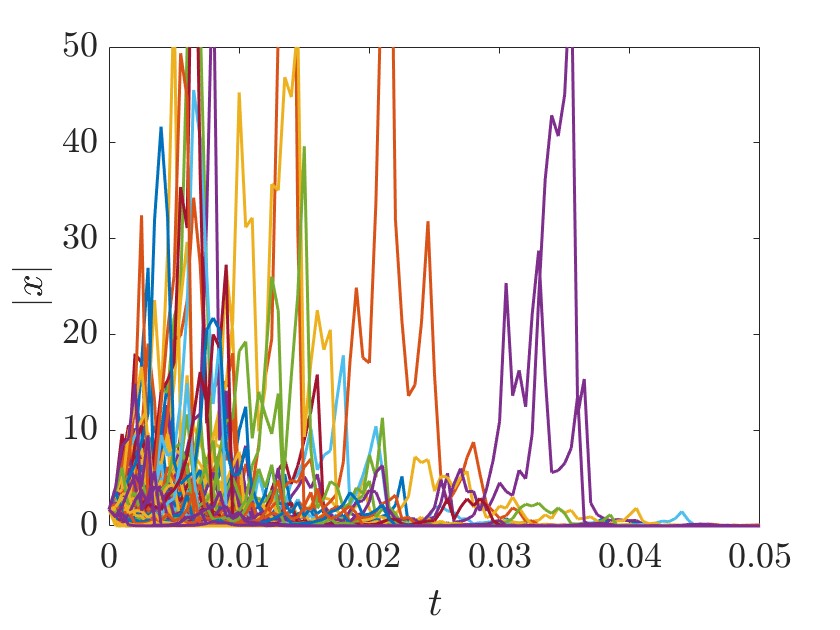}
		\label{figex4}}%
  \subfigure{
		\includegraphics[width=4cm]{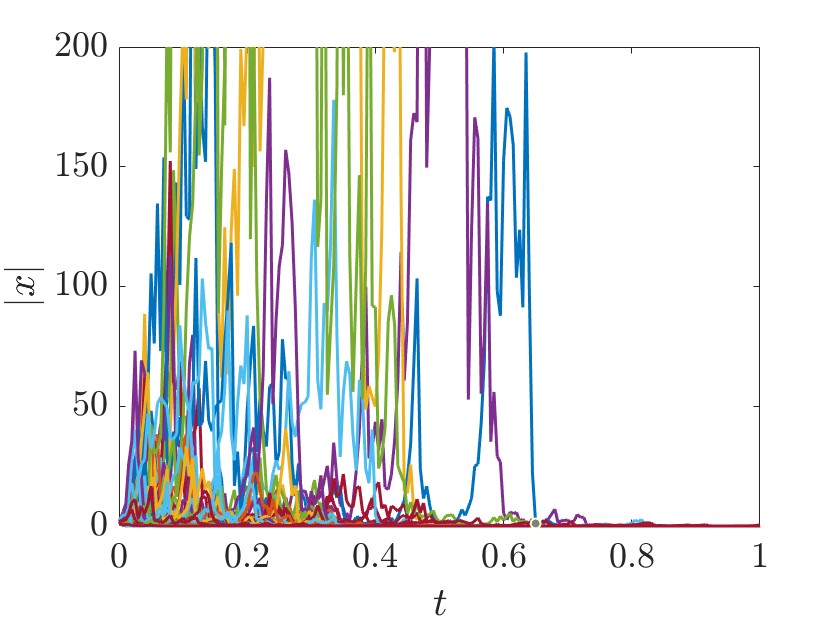}
		\label{fig3}}
	\centering
	\caption{(a)The dynamics of ${\rm log}\vert {\bm x}\vert $ change with $t$ for a group of SDEs correspondingly from the $G$-SDEs in Example~\ref{example1}.   Here, simulated are the 400 trials using the settings,  $\underline{\sigma}^2=3.5$, $\overline{\sigma}^2=4$.  (b)The dynamics of $\vert\bm x\vert$ change with $t$
		for a group of SDEs correspondingly from the $G$-SDEs in Example~\ref{exstability}.  Here, simulated are the 400 trials using the settings: $\underline{\sigma}^2=40$, $\overline{\sigma}^2=50$, $\sigma=10$, $\rho=10$, $\beta=8/3$, and $k=5$. (c)The dynamics of $\vert\bm x\vert$ change with $t$
		for a group of SDEs correspondingly from the $G$-SDEs in Example~\ref{example3}.  Here, simulated are the 400 trials using the settings: $\underline{\sigma}^2=40$ and $\overline{\sigma}^2=50$. }\label{fig1}
\end{figure}

In Figure \ref{fig1}, we show the numerical results, respectively, for Examples~\ref{example1}-\ref{example3}.

	%\bibliographystyle{aipsamp2}
%\bibliography{bibliography}	
	
	\bibliographystyle{siamplain}
\bibliography{bibliography}

	%\nocite{*}
	%\bibliography{aipsamp2}% Produces the bibliography via BibTeX.
	
\end{document}